\documentclass[ijoc,nonblindrev]{informs3} % current default for manuscript submission

\OneAndAHalfSpacedXII % current default line spacing
\usepackage{natbib}
 \bibpunct[, ]{(}{)}{,}{a}{}{,}%
 \def\bibfont{\small}%
 \def\bibsep{\smallskipamount}%
 \def\bibhang{24pt}%
 \def\BIBand{and}%

\TheoremsNumberedThrough     % Preferred (Theorem 1, Lemma 1, Theorem 2)
\EquationsNumberedThrough    % Default: (1), (2), ...
\usepackage{amsmath}
\usepackage{enumitem}
\usepackage{booktabs}
\usepackage{graphicx}
                 
\usepackage{bm}       
\usepackage{endnotes}
\let\footnote=\endnote

\usepackage[ruled,vlined]{algorithm2e}

\usepackage{listings}

\definecolor{codegreen}{rgb}{0,0.6,0}
\definecolor{codegray}{rgb}{0.5,0.5,0.5}
\definecolor{codepurple}{rgb}{0.58,0,0.82}
\definecolor{backcolour}{rgb}{0.95,0.95,0.92}

\usepackage{pdflscape}

\usepackage{multirow}

\newcommand{\pvtablestretch}{0.75}
\newcommand{\pveqnstretch}{1}

\newcommand{\field}[1]{\ensuremath{\mathbb{#1}}}
\newcommand{\sets}[1]{\ensuremath{\mathcal{#1}}}

\newcommand{\reals}{\ensuremath{\field{R}}} % real numbers
\newcommand{\integers}{\ensuremath{\field{Z}}} % integershttps://v2.overleaf.com/project/5b96bb3598d6c075a131c44a
\newcommand{\1}{\ensuremath{{\rm \mathbf e}}} % vector of all 1's
\newcommand{\I}[1]{\ensuremath{\mathbb{I}\left(#1\right)}} % indicator function
\newcommand{\subjectto}{\text{\rm subject to}} % subject to
\newcommand{\minimize}{\ensuremath{\mathop{\mathrm{minimize}}\limits}}
\newcommand{\maximize}{\ensuremath{\mathop{\mathrm{maximize}}\limits}}

\newcommand{\newpv}[1]{{\color{black}{#1}}}

\DeclareRobustCommand{\cpluspluslogo}{\hbox{C\hspace{-0.5ex}
                       \protect\raisebox{0.5ex}
                       {\protect\scalebox{0.67}{++}}}}

\def\Csharp{C\#}

\def\ROCPP{RO\cpluspluslogo}

\begin{document}
%%%%%%%%%%%%%%%%

%\setcounter{page}{1}
%\DOI{10.1287/xxxx.0000.0000}%

\RUNAUTHOR{Vayanos, Jin, and Elissaios}

\RUNTITLE{ROC++: Robust Optimization in C++}

% \TITLE{Bundling Information Goods of Decreasing Value}
\TITLE{RO\cpluspluslogo: Robust Optimization in \cpluspluslogo}

\ARTICLEAUTHORS{%
\AUTHOR{Phebe Vayanos,\textsuperscript{\textdagger} Qing Jin,\textsuperscript{\textdagger} and George Elissaios}
\AFF{\textsuperscript{\textdagger}University of Southern California, CAIS Center for Artificial Intelligence in Society \\ {\texttt{\{phebe.vayanos,qingjin\}@usc.edu}}}
%\AUTHOR{Author2}
%\AFF{Author2 affiliation, \EMAIL{}, \URL{}}
} % end of the block

\ABSTRACT{%
%\phantom{bla}\\
Over the last two decades, robust optimization techniques have emerged as a very popular means to address decision-making problems affected by uncertainty. Their success has been fueled by their attractive robustness and scalability properties, by ease of modeling, and by the limited assumptions they need about the uncertain parameters to yield meaningful solutions. Robust optimization techniques are available which can address both single- and multi-stage decision-making problems involving real-valued and/or binary decisions, and affected by both exogenous (decision-independent) and endogenous (decision-dependent) uncertain parameters. Many of these techniques apply to problems with either robust (worst-case) or stochastic (expectation) objectives and can thus be tailored to the risk preferences of the decision-maker. Robust optimization techniques rely on duality theory (potentially augmented with approximations) to transform a semi-infinite optimization problem to a finite program of benign complexity (the ``robust counterpart''). While writing down the model for a robust or stochastic optimization problem is usually a simple task, obtaining the robust counterpart requires expertise in robust optimization. To date, very few solutions are available that can facilitate the modeling and solution of such problems. This has been a major impediment to their being put to practical use. In this paper, we propose \ROCPP{}, a \cpluspluslogo~based platform for automatic robust optimization, applicable to a wide array of single- and multi-stage stochastic and robust problems with both exogenous and endogenous uncertain parameters. Our platform naturally extends existing off-the-shelf deterministic optimization platforms. We also propose the ROB file format that generalizes the LP file format to robust optimization. We showcase the modeling power of \ROCPP{} on several decision-making problems of practical interest. Our platform can help streamline the modeling and solution of stochastic and robust optimization problems for both researchers and practitioners. It comes with detailed documentation to facilitate its use and expansion. \ROCPP{} is freely distributed for academic use at \url{https://sites.google.com/usc.edu/robust-opt-cpp/}. % Enter your abstract
}%

%\KEYWORDS{deterministic inventory theory; infinite linear programming duality; 
%  existence of optimal policies; semi-Markov decision process; cyclic schedule}

\KEYWORDS{robust optimization, sequential decision-making under uncertainty, exogenous uncertainty, endogenous uncertainty, decision-dependent uncertainty, decision-dependent information discovery.}

\maketitle

\maketitle
%%%%%%%%%%%%%%%%%%%%%%%%%%%%%%%%%%%%%%%%%%%%%%%%%%%%%%%%%%%%%%%%%%%%%%

%\newpage

% Samples of sectioning (and labeling) in OPRE
% NOTE: (1) \section and \subsection do NOT end with a period
%       (2) \subsubsection and lower need end punctuation
%       (3) capitalization is as shown (title style).
%
%\section{Introduction.}\label{intro} %%1.
%\subsection{Duality and the Classical EOQ Problem.}\label{class-EOQ} %% 1.1.
%\subsection{Outline.}\label{outline1} %% 1.2.
%\subsubsection{Cyclic Schedules for the General Deterministic SMDP.}
%  \label{cyclic-schedules} %% 1.2.1
%\section{Problem Description.}\label{problemdescription} %% 2.

% Text of your paper here

%%%%%%%%%%%%%%%%%%%%%%%%%%%%%%%%%%%%%%%%%%%%%%%%%%%%%%%%%%%%%%%%%%%%%%%%%%%%%
%%%%%%%%%%%%%%%%%%%%%%%%%%%%%%%%%%%%%%%%%%%%%%%%%%%%%%%%%%%%%%%%%%%%%%%%%%%%%
%%%%%%%%%%%%%%%%%%%%%%%%%%%%%%%%%% INTRODUCTION %%%%%%%%%%%%%%%%%%%%%%%%%%%%% %%%%%%%%%%%%%%%%%%%%%%%%%%%%%%%%%%%%%%%%%%%%%%%%%%%%%%%%%%%%%%%%%%%%%%%%%%%%%
%%%%%%%%%%%%%%%%%%%%%%%%%%%%%%%%%%%%%%%%%%%%%%%%%%%%%%%%%%%%%%%%%%%%%%%%%%%%%

\section{Introduction}
\label{sec:introduction}

%%%%%%%%%%%%%%%%%%%%%%%%%%%%%%%%%%%%%%%%%%%%%%%%%%%%%%%%%%%%%%%%%%%%%%%%%%%%
%%%%%%%%%%%%%%%%%%%%%%%%%%%%%%%%%%%%%%%%%%%%%%%%%%%%%%%%%%%%%%%%%%%%%%%%%%%%

\subsection{Background \& Motivation}

Decision-making problems involving \emph{uncertain} or \emph{unknown} parameters are faced routinely by individuals, firms, policy-makers, and governments. 
%
%Decision-making problems affected by uncertainty are ubiquitous in the real-world where most decisions made by individuals, firms, policy-makers, and even governments involve at least some elements that are not known in advance. This is the case for example in inventory management where inventory levels (e.g., national stockpile of Personal Protective Equipment (PPE), ventilators, and other healthcare resources in preparation for a pandemic) must be decided before the demand for products becomes known (e.g., in advance of knowing the type of virus that will strike and the number of individuals that will get infected).% Other examples are in the context of financial planning where investments must be made before stock and bond price evolution is revealed, and in energy capacity expansion planning where plants must be built before availability of intermittent sources of electricity becomes known.
%
%
%Optimization under uncertainty is concerned with the modeling and solution of such decision-making problems involving parameters that are not perfectly known at the time when decisions are made. 
%
%Such parameters are referred to as \emph{uncertain} and 
Uncertain parameters may correspond to prediction errors, measurement errors, or implementation errors, see e.g., \cite{BenTal_Book}. Prediction errors arise when some of the data elements have not yet materialized at the time of decision-making and must thus be predicted/estimated (e.g., future prices of stocks, future demand, or future weather). Measurement errors arise when some of the data elements (e.g., characteristics of raw materials) cannot be precisely measured (due to e.g., limitations of the technological devices available). Implementation errors arise when some of the decisions may not be implemented exactly as planned/recommended by the optimization (due to e.g., physical constraints). 

If all decisions must be made \emph{before} the uncertain parameters are revealed, the decision-making problem is referred to as \emph{static} or \emph{single-stage}. In contrast, if the uncertain parameters are revealed sequentially over time and decisions are allowed to adapt to the history of observations, the decision-making problem is referred to as \emph{adaptive} or \emph{multi-stage}. 

In sequential decision-making problems, the time of revelation of the uncertain parameters may either be known a-priori or it may be part of the decision space. Uncertain parameters whose time of revelation is known in advance, being \emph{in}dependent of the decision-maker's actions, are referred to as \emph{exogenous}. Uncertain parameters whose time of revelation can be controlled by the decision-maker are referred to as \emph{endogenous}. This terminology was originally coined by~\cite{Jonsbraten_thesis}.

Examples of decision-making problems involving exogenous uncertain parameters are: financial portfolio optimization (see e.g., \cite{PortfolioSelection_Markowitz1952}), inventory and supply-chain management (see e.g., \cite{Scarf_Inventory}), vehicle routing (\cite{Bertsimas_1991_VRP}), unit commitment (see e.g., \cite{Takriti_1996}), and option pricing (see e.g., \cite{Kuhn_SwingOptions}). Examples of decision-making problems involving endogenous uncertain parameters are: R\&D project portfolio optimization (see e.g., \cite{Solak_RandDportfolios}), clinical trial planning (see e.g., \cite{Colvin_ClinicalTrials}), offshore oilfield exploration (see e.g., \cite{GoelGrossmanGasFields}), best box and Pandora's box problems (see e.g., \cite{Weitzman_1979}), and preference elicitation (see e.g., \cite{Vayanos_ActivePreferences}).

\subsection{Stochastic \& Robust Optimization}

Whether the decision-making problem is affected by exogenous and/or endogenous uncertain parameters, it is well known that ignoring uncertainty altogether when deciding on the actions to take usually results in suboptimal or even infeasible actions. To this end, researchers in stochastic and robust optimization have devised optimization-based models and solution approaches that explicitly capture the uncertain nature of these parameters. These frameworks model decisions as functions (\emph{decision rules}) of the history of observations, capturing the adaptive and non-anticipative nature of the decision-making process. 

Stochastic optimization, also known as stochastic programming, assumes that the distribution of the uncertain parameters is perfectly known, see e.g., \cite{StochasticProgrammingBook,SP_prekopa,Birge_Book}, and~\cite{SPshapiro_book}. This assumption is well justified in many situations. For example, this is the case if the distribution is stationary and can be well estimated from historical data. If the distribution of the uncertain parameters is discrete, the stochastic program admits a \emph{deterministic equivalent} that can be solved with off-the-shelf solvers potentially augmented with decomposition techniques, see e.g., \cite{benders1962}, or dedicated algorithms, see e.g., \cite{Rockafellar1991}. If the distribution of the uncertain parameters is continuous, the reformulation of the uncertain optimization problem may or not be computationally tractable since even evaluating the objective function usually requires computing a high-dimensional integral. If this problem is not computationally tractable, discretization approaches (such as the sample average approximation) may be employed. While discretization appears as a promising approach for
smaller problems, it may result in a combinatorial state explosion when applied to large and medium sized problems. Conversely, using only very few discretization points can result in solutions that are suboptimal or may even fail to be implementable in practice. Over the last two decades, stochastic programming techniques have been extended to address problems involving endogenous uncertain parameters, see e.g., \cite{GoelGrossmanGasFields,GoelGrossmann_BandB_DDU,GoelGrossman_ClassStochastic_DDU,GoelGrossman_NovelBB_GasFields,GuptaGrossman_SolutionStrategies,TarhanGrossmanGoel_Nonconvex_DDU,Colvin_ClinicalTrials,Colvin_TestingTasks,Colvin_Pharmaceutical}. We refer the reader to~\cite{StochasticProgrammingBook,SP_prekopa,Birge_Book}, and~\cite{SPshapiro_book} for in-depth reviews of the field of stochastic programming. 

Robust optimization does not necessitate knowledge of the distribution of the uncertain parameters. Rather than modeling uncertainty by means of distributions, it merely assumes that the uncertain parameters belong in a so-called \emph{uncertainty set}. The decision-maker then seeks to be immunized against all possible realizations of the uncertain parameters in this set. The robust optimization paradigm gained significant traction starting in the late 1990s and early 2000s following the works of \cite{RSols_UncertainLinear_Programs_BT,BenTalNemirovski_RCO,RobustLP_UncertainData,AdjustableRSols_uncertainLP}, and~\cite{bertsimas2003robust,Price_Robustness,Tractable_R_SOCP}, among others. Over the last two decades, research on robust optimization has burgeoned, fueled by the limited assumptions it needs about the uncertain parameters to yield meaningful solutions, by its attractive robustness and scalability properties, and by ease of modelling, see e.g., \cite{TheoryApplicationRO,GYdH15:practical_guide}.%Gabrel2014

Robust optimization techniques are available which can address both single- and multi-stage decision-making problems involving real-valued and/or binary decisions, and affected by exogenous and/or endogenous uncertain parameters. Robust optimization techniques for single-stage robust optimization rely on duality theory to transform a semi-infinite optimization problem to an equivalent finite program of benign complexity (the ``robust counterpart'') that is solvable with off-the-shelf solvers, see e.g.,~\cite{BenTal_Book}. In the multi-stage setting, the dualization step is usually preceded by an approximation step that transforms the multi-stage problem to a single-stage robust program. The idea is to restrict the space of the decisions to a subset of benign complexity based either on a decision rule approximation or a finite adaptability approximation. The decision rule approximation consists in restricting the adjustable decisions to those presenting e.g., linear, piecewise linear, or polynomial dependence on the uncertain parameters, see e.g.,~\cite{AdjustableRSols_uncertainLP,kuhn_primal_dual_rules,Bertsimas_polynomial_policies,ConstraintSampling_VKR,Angelos_Liftings}. The finite adaptability approximation consists in selecting a finite number of candidate strategies today and selecting the best of those strategies in an adaptive fashion once the uncertain parameters are revealed, see e.g., \cite{Caramanis_FiniteAdaptability,Hanasusanto2015}. The decision rule and finite adaptability approximations have been extended to the endogenous uncertainty setting, see~\cite{DDI_VKB,vayanos_ROInfoDiscovery}. While writing down the model for a robust optimization problem is usually a simple task (akin to formulating a deterministic optimization problem), obtaining the robust counterpart is typically tedious and requires expertise in robust optimization, see~\cite{BenTal_Book}.

Robust optimization techniques have been extended to address certain classes of stochastic programming problems involving continuously distributed uncertain parameters and affected by both exogenous uncertainty (see \cite{kuhn_primal_dual_rules,Bodur2018}) and endogenous uncertainty (see~\cite{DDI_VKB}).

Robust optimization techniques have been used successfully to address single-stage problems in inventory management (\cite{AJD16:inventory}), network optimization (\cite{bertsimas2003robust}), product pricing (\cite{Adida_DP_MP,Thiele_MPP_JRPM}), portfolio optimization (\cite{Robust_Portfolio_Management,RobustPortfolioSelection}), and healthcare (\cite{gupta2017maximizing,BTV_kidneys,ChanDefibrillator2018}). They have also been used to successfully tackle sequential problems in energy (\cite{ZWWG13:multistage_robust_unit_commitment,Jiang2014}), inventory and supply-chain management (\cite{BenTal_RSFC,mamani2016closed}), network optimization (\cite{Atamturk_NetworkDesign}), preference elicitation (\cite{Vayanos_ActivePreferences}), vehicle routing (\cite{Gounaris_RobustVehicleRouting}), process scheduling (\cite{Lappas2016}), and R\&D project portfolio optimization (\cite{vayanos_ROInfoDiscovery}).

In spite of its success at addressing a diverse pool of problems in the literature, to date, very few platforms are available that can facilitate the modeling and solution of robust optimization problems, and those available can only tackle limited classes of robust problems. At the same time, and as mentioned above, reformulating such problems in a way that they can be solved by off-the-shelf solvers requires expertise. This is particularly true in the case of multi-stage problems and of problems affected by endogenous uncertainty. 

In this paper, we fill this gap by proposing \ROCPP, a \cpluspluslogo~based platform for modeling, approximating, automatically reformulating, and solving general classes of robust optimization problems. Our platform provides several modeling objects (decision variables, uncertain parameters, constraints, optimization problems) and overloaded operators to allow for ease of modeling using a syntax similar to that of state-of-the-art \emph{deterministic optimization} solvers like CPLEX\footnote{\url{https://www.ibm.com/analytics/cplex-optimizer}} or Gurobi.\footnote{\url{https://www.gurobi.com}} While our platform is not exhaustive, it provides a framework that is easy to update and expand and lays the foundation for more development to help facilitate research in, and real-life applications of, robust optimization.

%Over the past decade, two key impediments have been informally recognized as central in inhibiting the widespread industrial use of stochastic programming. First, modeling systems for mathematical programming have only recently begun to incor- porate extensions for specifying stochastic programs. Without an integrated and acces- sible modeling capability, practitioners are forced to implement custom techniques for specifying stochastic programs. Second, stochastic programs are often extremely dif- ficult to solve—especially in contrast to their deterministic counterparts. There exists no analog to CPLEX [12], GUROBI [25], or XpressMP [60] for stochastic program- ming, principally because the algorithmic technology is still under active investigation and development, particularly in the multi-stage, non-linear, and mixed-integer cases.

%%%%%%%%%%%%%%%%%%%%%%%%%%%%%%%%%%%%%%%%%%%%%%%%%%%%%%%%%%%%%%%%%%%%%%%%%%%%%%%%%%%%%%
%%%%%%%%%%%%%%%%%%%%%%%%%%%%%%%%%%%%%%%%%%%%%%%%%%%%%%%%%%%%%%%%%%%%%%%%%%%%%%%%%%%%%%

\subsection{Related Literature}

% Xpress.\footnote{See \url{https://www.fico.com/en/products/fico-xpress-solver}.} 
% Knitro,\footnote{See \url{https://www.artelys.com/solvers/knitro/}.}

\paragraph{Tools for Modelling and Solving Deterministic Optimization Problems.} There exist many commercial and open-source tools for modeling and solving deterministic optimization problems. On the commercial front, the most popular solvers for conic (integer) optimization are~Gurobi,\footnote{See \url{https://www.gurobi.com}.} IBM CPLEX Optimizer,\footnote{See \url{https://www.ibm.com/analytics/cplex-optimizer}.} and Mosek.\footnote{See \url{https://www.mosek.com}.} These solvers provide interfaces for C/\cpluspluslogo, Python, and other commonly used high-level languages. Several tools such as AMPL,\footnote{See \url{https://ampl.com}} GAMS,\footnote{See \url{https://www.gams.com}.} and AIMMS\footnote{See \url{https://www.aimms.com}.} are based on dedicated modelling languages. They provide APIs for \cpluspluslogo, \Csharp, Java, and Python. They can also connect to commercial or open-source solvers. Finally, several commercial vendors provide modeling capabilities combined with built-in solvers, see e.g., Lindo Systems Inc.,\footnote{See \url{https://www.lindo.com}.} FrontlineSolvers,\footnote{See \url{https://www.solver.com/}.} and Maximal.\footnote{See \url{http://www.maximalsoftware.com/}.} On the open source side, the most popular solvers are GLPK,\footnote{See \url{https://www.gnu.org/software/glpk/}.} and Cbc\footnote{See \url{https://github.com/coin-or/Cbc}} and Clp\footnote{See \url{https://github.com/coin-or/Clp}} from COIN-OR.\footnote{See \url{https://www.coin-or.org}} Commercial and open-source solvers can also be accessed from several open-source modeling languages for mathematical optimization that are embedded in popular algebraic modeling languages. These include JuMP which is embedded in Julia, see \cite{JuMP}, and Yalmip and CVX which are embedded in MATLAB, see \cite{YALMIP} and \cite{cvx,gb08}, respectively.

\paragraph{Tools for Modelling and Solving Stochastic Optimization Problems.} Several of the commercial vendors also provide modeling capabilities for stochastic programming, see e.g., Lindo Systems Inc.,\footnote{See \url{https://www.lindo.com/index.php/products/lingo-and-optimization-modeling}.} FrontlineSolvers,\footnote{See \url{https://www.solver.com/risk-solver-stochastic-libraries}.} Maximal,\footnote{See \url{http://www.maximalsoftware.com/maximal/news/stochastic.html}.} GAMS,\footnote{See \url{https://www.gams.com/latest/docs/UG_EMP_SP.html}.} AMPL (see \cite{SAMPL}), and AIMMS.\footnote{See \url{https://www.aimms.com}.} On the other hand, there are only two open-source platforms that we are aware of that provide such capabilities. The first one is FLOP\cpluspluslogo,\footnote{See \url{https://projects.coin-or.org/FlopC++}.} which is part of COIN-OR. It provides an algebraic modeling environment in \cpluspluslogo\ that is similar to languages such as GAMS and AMPL. The second one is PySP\footnote{See \url{https://pyomo.readthedocs.io/en/stable/modeling_extensions/pysp.html}.} which is based on the Python high-level programming language, see \cite{pysp2012}. To express a stochastic program in PySP, the user specifies both the deterministic base model and the scenario tree model in the Pyomo open-source algebraic modeling language, see~\cite{pyomo2009,pyomo2011,pyomo2012}. All the aforementioned tools assume that the distribution of the uncertain parameters in the optimization problem is discrete or provide mechanisms for generating samples from a continuous distribution to feed into the model.

\paragraph{Tools for Modelling and Solving Robust Optimization Problems.} Our robust optimization platform \ROCPP{} most closely relates to several open-source tools released in recent years. All of these tools present a similar structure: they provide a modeling platform combined with an approximation/reformulation toolkit that can automatically obtain the robust counterpart, which is then solved using existing open-source and/or commercial solvers. The majority of these platforms is based on the MATLAB modeling language. One tool builds upon YALMIP, see~\cite{lofberg2012}, and provides support for single-stage problems with exogenous uncertainty. A notable advantage of YALMIP is that the robust counterpart output by the platform can be solved using any one of a huge variety of open-source or commercial solvers. Other platforms, like ROME and RSOME are entirely motivated by the (stochastic) robust optimization modeling paradigm, see~\cite{ROME} and \cite{RSOME}, and provide support for both single- and multi-stage (distributionally) robust optimization problems affected by exogenous uncertain parameters. The robust counterparts output by ROME and RSOME can be solved with CPLEX, Mosek, and SDPT3.\footnote{See \url{http://www.math.cmu.edu/~reha/sdpt3.html}.}  Recently, JuMPeR has been proposed as an add-on to JuMP. It can cater for single-stage problems with exogenous uncertain parameters. JuMPeR can be connected to a large variety of open-source and commercial solvers. On the commercial front, AIMMS is currently equipped with an add-on that can be used to model and automatically reformulate robust optimization problems. It can tackle both single- and multi-stage problems with exogenous uncertainty. A CPLEX license is needed to operate this add-on. To the best of our knowledge, none of the available platforms can address problems involving endogenous uncertain parameters. None of them can tackle problems presenting binary adaptive variables. Finally, none of these platforms can be used from \cpluspluslogo.

\paragraph{File Formats for Specifying Optimization Problems.} To facilitate the sharing and storing of optimization problems, dedicated file formats have been proposed. The two most popular file formats for deterministic mathematical programming problems are the MPS and LP formats. MPS is an older format established on mainframe systems. It is not very intuitive to use as it is setup as if you were using punch cards. In contrast, the LP format is a lot more interpretable: it captures problems in a way similar to how it is modelled on paper. The SMPS file format is the most popular format for storing stochastic programs and mirrors the role MPS plays in the deterministic setting, see \cite{BirgDempGassGunnKingWall87,Gassmann_SLP}. To the best of our knowledge, no format exists in the literature for storing and sharing robust optimization problems.

\subsection{Contributions}

We now summarize our main contributions and the key advantages of our platform:
\begin{enumerate}[label=\textit{(\alph*)}]
\item We propose \ROCPP{}, the first \cpluspluslogo\ based platform for modelling, automatically reformulating, and solving robust optimization problems. Our platform is the first capable of addressing both single- and multi-stage problems involving exogenous and/or endogenous uncertain parameters and real- and/or binary-valued adaptive variables. It can also be used to address certain classes of single- or multi-stage stochastic programs whose distribution is continuous and supported on a compact set. Our reformulations are (mixed-integer) linear or second-order cone optimization problems and thus any solver that can tackle such problems can be used to solve the robust counterparts output by our platform. We provide an interface to the commercial solver Gurobi. Our platform can easily be extended to support other solvers. We illustrate the flexibility and ease of use of our platform on several stylized problems. 
\item We propose the ROB file format, the first file format for storing and sharing general robust optimization problems. Our format builds upon the LP file format and is thus interpretable and easy to use.
\item Our modeling language is similar to the one provided for the deterministic case by solvers such as CPLEX or Gurobi: it is easy to use for anyone familiar with these.
\item Our platform comes with detailed documentation (created with Doxygen\footnote{See \url{https://www.doxygen.nl/index.html}.}) to facilitate its use and expansion. Our framework is open-source for educational, research, and non-profit purposes. The source code, installation instructions, and dependencies of \ROCPP{} are available at \url{https://sites.google.com/usc.edu/robust-opt-cpp/}.
\end{enumerate}

%%%%%%%%%%%%%%%%%%%%%%%%%%%%%%%%%%%%%%%%%%%%%%%%%%%%%%%%%%%%%%%%%%%%%%%%%%%%%%%%%%%%%%
%%%%%%%%%%%%%%%%%%%%%%%%%%%%%%%%%%%%%%%%%%%%%%%%%%%%%%%%%%%%%%%%%%%%%%%%%%%%%%%%%%%%%%

\subsection{Organization of the Paper \& Notation}

The remainder of this paper is organized as follows. Section~\ref{sec:modelling_background} describes the broad class of problems to which \ROCPP{} applies. Section~\ref{sec:uncertainty_set} presents our model of uncertainty. Section~\ref{sec:approximation_schemes} lists the approximation schemes that are provided by \ROCPP. Sample models created and solved using \ROCPP{} are provided in Section~\ref{sec:numerical_results}. Section~\ref{sec:file_format} introduces the ROB file format. Section~\ref{sec:extensions} presents extensions to the core model that can also be tackled by \ROCPP.

\paragraph{Notation.} Throughout this paper, vectors (matrices) are denoted by boldface lowercase (uppercase) letters. The $k$th element of a vector ${\bm x} \in \reals^n$ ($k \leq n$) is denoted by ${\bm x}_k$. Scalars are denoted by lowercase or upper case letters, e.g., $\alpha$ or $N$. We let $\mathcal L_{k}^n$ denote the space of all functions from $\reals^k$ to $\reals^n$. Accordingly, we denote by $\mathcal B_{k}^n$ the space of all functions from $\reals^k$ to $\{0,1\}^n$. Given two vectors of equal length, ${\bm x}$, ${\bm y} \in \mathbb R^n$, we let ${\bm x} \circ {\bm y}$ denote the Hadamard product of the vectors, i.e., their element-wise product. 

Throughout the paper, we denote the uncertain parameters by ${\bm \xi} \in \reals^k$. We consider two settings: a robust setting and a stochastic setting. In the robust setting, we assume that the decision-maker wishes to be immunized against realizations of ${\bm \xi}$ in the uncertainty set $\Xi$. In the stochastic setting, we assume that the distribution $\mathbb P$ of the uncertain parameters is fully known. In this case, we let $\Xi$ denote its support and we let $\mathbb E(\cdot)$ denote the expectation operator with respect to $\mathbb P$.

%%%%%%%%%%%%%%%%%%%%%%%%%%%%%%%%%%%%%%%%%%%%%%%%%%%%%%%%%%%%%%%%%%%%%%%%%%%%%%%%%%%%%%
%%%%%%%%%%%%%%%%%%%%%%%%%%%%%%%%%%%%%%%%%%%%%%%%%%%%%%%%%%%%%%%%%%%%%%%%%%%%%%%%%%%%%%
%%%%%%%%%%%%%%%%%%%%%%%%%%%%%% MODELLING WITH RO %%%%%%%%%%%%%%%%%%%%%%%%%%%%%%%%%%%%% %%%%%%%%%%%%%%%%%%%%%%%%%%%%%%%%%%%%%%%%%%%%%%%%%%%%%%%%%%%%%%%%%%%%%%%%%%%%%%%%%%%%%%
%%%%%%%%%%%%%%%%%%%%%%%%%%%%%%%%%%%%%%%%%%%%%%%%%%%%%%%%%%%%%%%%%%%%%%%%%%%%%%%%%%%%%%

\section{Modelling Decision-Making Problems Affected by Uncertainty}
\label{sec:modelling_background}

We present the main class of problems that are supported by \ROCPP. Our platform can handle general multi-stage decision problems affected by both exogenous and endogenous uncertainty over the finite planning horizon $\mathcal T: = \{1,\ldots,T\}$ (if $T=1$, then we recover single-stage decision-problems). The elements of the vector of uncertain parameters ${\bm \xi}$ are revealed sequentially over time. However, the sequence of their revelation need not be predetermined (exogenous). Instead, the time of information discovery can be controlled by the decision-maker via the binary \emph{measurement decisions} ${\bm w}$. The aim is to find the sequences of real-valued decisions ${\bm y} := ({\bm y}_1,\ldots,{\bm y}_T)$, binary decisions ${\bm z}:=({\bm z}_1,\ldots,{\bm z}_T)$, and measurement decisions ${\bm w}:=({\bm w}_1,\ldots,{\bm w}_T)$ that minimize a given (uncertain) cost function either in expectation (stochastic setting) or in the worst-case (robust setting). These decisions are constrained by a set of inequalities which are required to be obeyed robustly, i.e., for any realization of the parameters~${\bm \xi}$ in the set $\Xi$. A salient feature of our platform is that the decision variables are explicitly modelled as functions, or decision rules, of the history of observations. This feature is critical as it captures the ability of the decision-maker to adjust their decisions based on the realizations of the observed uncertain parameters.

Decision-making problems of the type described
here can be formulated as 
\begin{equation}\renewcommand{\arraystretch}{\pveqnstretch}
    \begin{array}{cl}
       \minimize  & \quad \displaystyle \mathbb F \left[ \;\; \sum_{t \in \sets T} {\bm c}_t^\top {\bm y}_t({\bm \xi}) + {\bm d}_t({\bm \xi}) {\bm z}_t({\bm \xi}) + {\bm f}_t({\bm \xi}) {\bm w}_t({\bm \xi}) \right]  \\
       \subjectto  & \quad {\bm y}_t \in \sets L_k^{n_t}, \; {\bm z}_t \in \sets B_k^{\ell_t}, \; {\bm w}_t \in \sets B_k^k \quad \forall t \in \sets T \\
       & \quad
       \!\!\!\!\left. \begin{array}{l}
       \displaystyle \sum_{\tau=1}^t {\bm A}_{t \tau}{\bm y}_\tau({\bm \xi}) + {\bm B}_{t \tau}({\bm \xi}){\bm z}_\tau({\bm \xi}) + {\bm C}_{t \tau}({\bm \xi}) {\bm w}_\tau({\bm \xi}) \; \leq \; h_t ({\bm \xi}) \\
       {\bm w}_t({\bm \xi}) \in \sets W_t \\
       {\bm w}_t({\bm \xi}) \geq {\bm w}_{t-1}({\bm \xi}) 
       \end{array} \quad \right\} \quad \forall {\bm \xi} \in \Xi, t \in \sets T \\
       & \quad
       \!\!\!\!\left. \begin{array}{l}
        {\bm y}_t({\bm \xi}) = {\bm y}_t({\bm \xi}') \\
        {\bm z}_t({\bm \xi}) = {\bm z}_t({\bm \xi}') \\
        {\bm w}_t({\bm \xi}) = {\bm y}_t({\bm \xi}')
        \end{array} \quad \quad  \right\}  \quad \forall t \in \sets T,\; \forall {\bm \xi},\; {\bm \xi}' \in \Xi : {\bm w}_{t-1}({\bm \xi}) \circ {\bm \xi} = {\bm w}_{t-1}({\bm \xi}') \circ {\bm \xi}',
    \end{array}
\label{eq:main_problem}
\end{equation}
where $\mathbb F$ is a functional that maps the uncertain overall costs (across all possible realizations of ${\bm \xi}$) to a real number. In the robust setting, the functional $\mathbb F(\cdot)$ computes the worst-case (maximum) over ${\bm \xi} \in \Xi$. In the stochastic setting, it computes the expectation of the cost function under the distribution of the uncertain parameters (assumed to be known).

In this problem, ${\bm y}_t({\bm \xi}) \in \reals^{n_t}$ (resp.\ ${\bm z}_t({\bm \xi}) \in \{0,1\}^{\ell_t}$) represent real-valued (resp.\ binary) decisions that are selected at the beginning of time period $t$. The variables ${\bm w}_t({\bm \xi}) \in \{0,1\}^k$ are binary measurement decisions that are made at time $t$ and that summarize the information base for time $t+1$. Specifically, ${\bm w}_{t}({\bm \xi})$ is a binary vector that has the same dimension as the vector of uncertain parameters; its $i$th element, ${\bm w}_{t,i}$, equals one if and only if the $i$th uncertain parameter, ${\bm \xi}_i$, has been observed at some time $\tau \in \{0,\ldots,t\}$, in which case it is included in the information basis for time $t+1$. Without loss of generality, we assume that ${\bm w}_0({\bm \xi})={\bm 0}$ for all ${\bm \xi} \in \Xi$ so that no uncertain parameter is known at the beginning of the planning horizon. The costs associated with the variables ${\bm y}_t({\bm \xi})$, ${\bm z}_t({\bm \xi})$, and ${\bm w}_t({\bm \xi})$ are ${\bm c}_t \in \reals^{n_t}$, ${\bm d}_t({\bm \xi}) \in \reals^{\ell_t}$, and ${\bm f}_t({\bm \xi}) \in \reals^k$, respectively. In particular, ${\bm f}_{t,i}({\bm \xi})$ represents the cost of including ${\bm \xi}_i$ in the information basis at time $t+1$. Without much loss of generality, we assume that the costs ${\bm d}_t({\bm \xi})$ and ${\bm f}_t({\bm \xi})$ are all linear in ${\bm \xi}$.

The first set of constraints in the formulation above defines the decision variables of the problem and ensures that decisions are modelled as functions of the uncertain parameters. The second set of constraints, which involve the matrices ${\bm A}_{t\tau} \in \reals^{m_t \times n_\tau}$, ${\bm B}_{t\tau}({\bm \xi}) \in \reals^{m_t \times \ell_\tau}$, and ${\bm C}_{t\tau}({\bm \xi}) \in \reals^{m_t \times k}$ are the problem constraints. The set $\mathcal W_t$ in the third constraint may model requirements stipulating, for example, that a specific uncertain parameter can only be observed after a certain stage. If the $i$th uncertain parameter is exogenous (i.e., if its time of information discovery is \emph{not} decision-dependent) and if its time of revelation is $t$, it suffices to set ${\bm w}_{\tau,i}({\bm \xi})=0$, if $\tau <t$; $=1$, else, and ${\bm f}_{\tau,i}({\bm \xi})=0$ for all ${\bm \xi} \in \Xi$ and $\tau \in \sets T$. The fourth set of constraints is an information monotonicity constraint: it stipulates that information that has been observed cannot be forgotten. The last three sets of constraints are decision-dependent non-anticipativity constraints: they stipulate that the adaptive decision-variables must be constant in those uncertain parameters that have not been observed at the time when the decision is made. Without much loss of generality, we assume that the matrices ${\bm B}_{t\tau}({\bm \xi}) \in \reals^{m_t \times \ell_\tau}$ and ${\bm C}_{t\tau}({\bm \xi}) \in \reals^{m_t \times k}$ are both linear in ${\bm \xi}$.

%%%%%%%%%%%%%%%%%%%%%%%%%%%%%%%%%%%%%%%%%%%%%%%%%%%%%%%%%%%%%%%%%%%%%%%%%%%%%%%%%%%%%%
%%%%%%%%%%%%%%%%%%%%%%%%%%%%%%%%%%%%%%%%%%%%%%%%%%%%%%%%%%%%%%%%%%%%%%%%%%%%%%%%%%%%%%
%%%%%%%%%%%%%%%%%%%%%%%%%%%% UNCERTAINTY MODELLING %%%%%%%%%%%%%%%%%%%%%%%%%%%%%%%%%%% %%%%%%%%%%%%%%%%%%%%%%%%%%%%%%%%%%%%%%%%%%%%%%%%%%%%%%%%%%%%%%%%%%%%%%%%%%%%%%%%%%%%%%
%%%%%%%%%%%%%%%%%%%%%%%%%%%%%%%%%%%%%%%%%%%%%%%%%%%%%%%%%%%%%%%%%%%%%%%%%%%%%%%%%%%%%%

\section{Modelling Uncertainty}
\label{sec:uncertainty_set}

We now discuss our model for the set $\Xi$. Throughout this paper and in our platform, we assume that $\Xi$ is compact and admits a conic representation, i.e., it is expressible as
\begin{equation}
\Xi \; := \left\{ {\bm \xi} \in \reals^k \; : \; \exists {\bm \zeta}^s \in \reals^{k_s}, \; s=1,\ldots, S \; : \; {\bm P}^s {\bm \xi} + {\bm Q}^s {\bm \zeta}^s + {\bm q}^s \in \sets K^s, \; s=1,\ldots, S \right\}
\label{eq:uncertainty_set}
\end{equation}
for some matrices ${\bm P}^s \in \reals^{r_s \times k}$ and ${\bm Q}^s \in \reals^{r_s \times k_s}$, and vector ${\bm q}^s \in \reals^{r_s}$, $s=1\ldots,S$, where $\sets K^s$, $s=1\ldots,S$, are closed convex pointed cones in $\reals^{r_s}$. Finally, we assume that the representation above is strictly feasible (unless the cones involved in the representation are polyhedral, in which case this assumption can be relaxed). In our platform, we focus on the cases where the cones $\sets K^s$ are either polyhedral, i.e., $\sets K^s = \reals_+^{r_s}$, or Lorentz cones, i.e., 
$
\sets K^s = \left\{ {\bm u} \in \reals^{r_s} \; :\; \sqrt{ {\bm u}_1^2 + \cdots + {\bm u}_{r_s-1}^2 } \leq {\bm u}_{r_s} \right\}.
$ 

Uncertainty sets of the form~\eqref{eq:uncertainty_set} arise naturally from statistics or from knowledge of the distribution of the uncertain parameters. In the stochastic setting, the uncertainty set $\Xi$ can be constructed as the support of the distribution of the uncertain parameters, see e.g., \cite{kuhn_primal_dual_rules}. More often, it is constructed in a data-driven fashion to guarantee that constraints are satisfied with high probability, see e.g., \cite{Bertsimas2018}. More generally, disciplined methods for constructing uncertainty sets from ``random'' uncertainty exist, see e.g., \cite{bandi2012tractable,BenTal_Book}. We now discuss several uncertainty sets from the literature that can be modelled in the form~\eqref{eq:uncertainty_set}.
\begin{example}[Budget Uncertainty Sets]
Uncertainty sets of the form~\eqref{eq:uncertainty_set} can be used to model 1-norm and $\infty$-norm uncertainty sets with budget of uncertainty $\Gamma$, given by $\{ {\bm \xi} \in \reals^k : \|{\bm \xi}\|_1 \leq \Gamma \}$ and $\{ {\bm \xi} \in \reals^k : \|{\bm \xi}\|_\infty \leq \Gamma \}$, respectively. More generally, they can be used to impose budget constraints at various levels of a given hierarchy. For example, they can be used to model uncertainty sets of the form
$$
\left\{ {\bm \xi} \in \reals^k \; : \; \sum_{i \in \sets H_h} |{\bm \xi}_i| \leq \Gamma_h {\bm \xi}_h \quad \forall h=1,\ldots, H \right\},
$$
where the sets $\sets H_h \subseteq \{1,\ldots,k\}$ collect the indices of all uncertain parameters in the $h$th level of the hierarchy and $\Gamma_h \in \reals_+$ is the budget of uncertainty for hierarchy $h$, see e.g., \cite{DSL2019}.
\end{example}

\begin{example}[Ellipsoidal Uncertainty Sets]
Uncertainty sets of the form~\eqref{eq:uncertainty_set} capture as special cases ellipsoidal uncertainty sets, which arise for example as confidence regions from Gaussian distributions. These are expressible as
$$
\left\{ {\bm \xi} \in \reals^k \; : \; ({\bm \xi}-\overline {\bm \xi})^\top {\bm P}^{-1} ({\bm \xi}-\overline {\bm \xi}) \; \leq \; 1 \right\},
$$
for some matrix ${\bm P} \in \mathbb S_+^{k}$ and vector $\overline {\bm \xi} \in \reals^k$, see e.g., \cite{BenTal_Book}.%symmetric positive definite
\end{example}

\begin{example}[Central Limit Theorem Uncertainty Sets]
Sets of the form~\eqref{eq:uncertainty_set} can be used to model Central Limit Theorem based uncertainty sets. These sets arise for example as confidence regions for large numbers of i.i.d.\ uncertain parameters and are expressible as
$$
\left\{ {\bm \xi} \in \reals^k : \left| \sum_{i=1}^k {\bm \xi}_k - \mu k \right| \leq \Gamma \sigma \sqrt{k}  \right\},
$$ 
where $\mu$ and $\sigma$ are the mean and standard deviation of the i.i.d.\ parameters ${\bm \xi}_i$, $i=1,\ldots,k$, see~\cite{bandi2012tractable,BTV_kidneys}.
\end{example}

\begin{example}[Uncertainty Sets based on Factor Models]
Sets of the form~\eqref{eq:uncertainty_set} capture as special cases uncertainty sets based on factor models that are popular in finance and economics. These are expressible in the form
$$
\left\{ {\bm \xi} \in \reals^k \; : \; \exists{\bm \zeta} \in \reals^{\kappa} : {\bm \xi} = {\bm \Phi} {\bm \zeta} + {\bm \phi}, \; \| {\bm \zeta}\|_2 \leq 1 \right\},
$$
for some vector ${\bm \phi} \in \reals^{k}$ and matrix ${\bm \Phi} \in \reals^{k \times \kappa}$.
\end{example}

%%%%%%%%%%%%%%%%%%%%%%%%%%%%%%%%%%%%%%%%%%%%%%%%%%%%%%%%%%%%%%%%%%%%%%%%%%%%%%%%%%%%%%
%%%%%%%%%%%%%%%%%%%%%%%%%%%%%%%%%%%%%%%%%%%%%%%%%%%%%%%%%%%%%%%%%%%%%%%%%%%%%%%%%%%%%%
%%%%%%%%%%%%%%%%%%%%%%%%%%%% APPROXIMATION SCHEMES %%%%%%%%%%%%%%%%%%%%%%%%%%%%%%%%%%% %%%%%%%%%%%%%%%%%%%%%%%%%%%%%%%%%%%%%%%%%%%%%%%%%%%%%%%%%%%%%%%%%%%%%%%%%%%%%%%%%%%%%%
%%%%%%%%%%%%%%%%%%%%%%%%%%%%%%%%%%%%%%%%%%%%%%%%%%%%%%%%%%%%%%%%%%%%%%%%%%%%%%%%%%%%%%

\section{Interpretable Decision Rules \& Contingency Planning}
\label{sec:approximation_schemes}

In formulation~\eqref{eq:main_problem}, the recourse decisions are very hard to interpret since decisions are modelled as (potentially very complicated) functions of the history of observations. The functional form of the decisions combined with the infinite number of constraints involved in problem~\eqref{eq:main_problem} also imply that this problem cannot be solved directly. This has motivated researchers in the fields of stochastic and robust optimization to propose several tractable approximation schemes capable of bringing problem~\eqref{eq:main_problem} to a form amenable to solution by off-the-shelf solvers. Broadly speaking, these approximation schemes fall in two categories: interpretable decision rule approximations which restrict the functional form of the recourse decisions; and finite adaptability approximation schemes which yield a moderate number of contingency plans that are candidates to be implemented in the future. These approximations have the benefit of improving the tractability properties of the problem and of rendering decisions more \emph{interpretable}, which is a highly desirable property of any automated decision support system.

\begin{figure}[ht!]
    \centering 
    \includegraphics[width=0.55\textwidth]{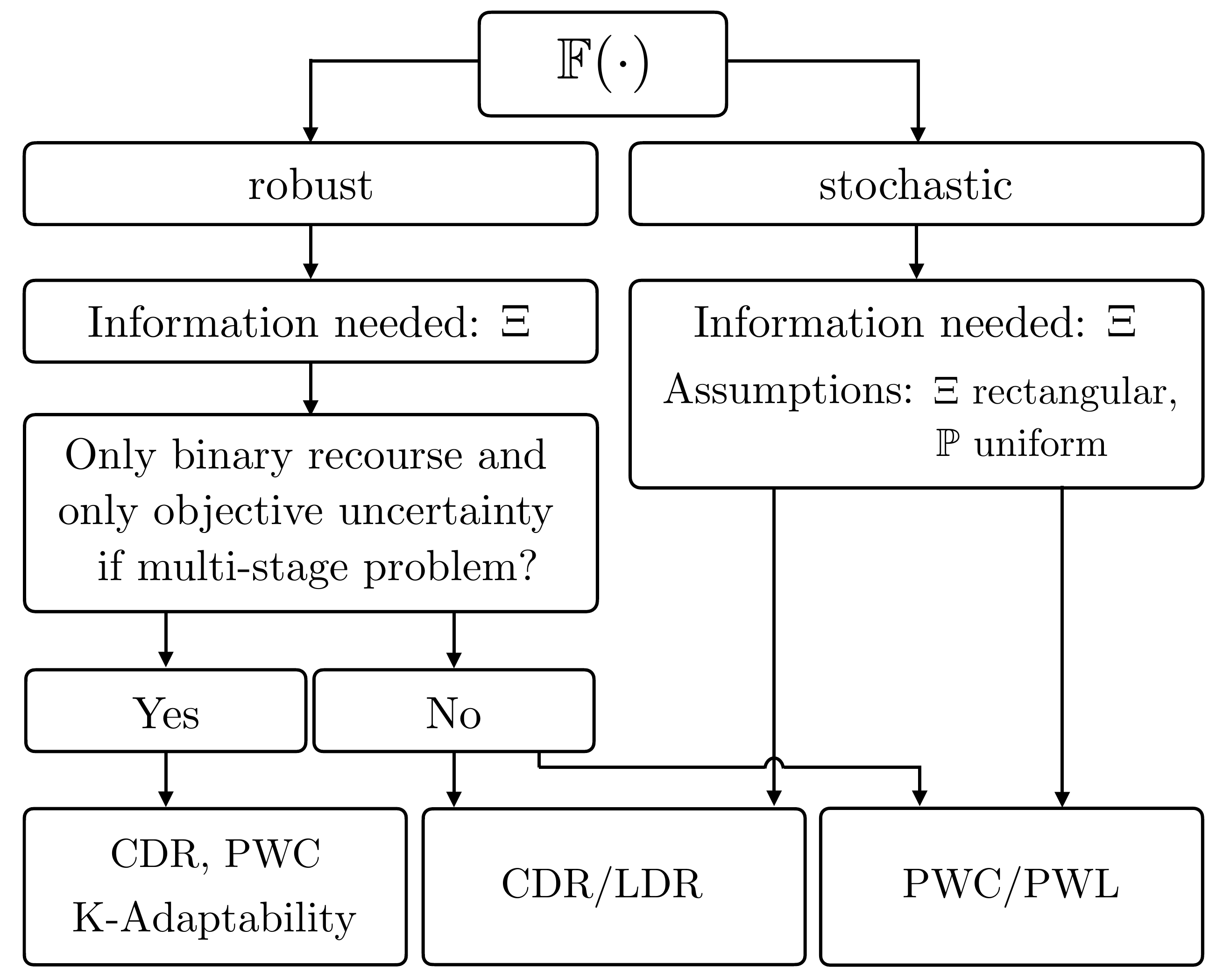}
    \caption{Decision tree to help guide the choice of approximation scheme for multi-stage problems in \ROCPP{}.}
    \label{fig:decision_tree}
\end{figure}

We now describe the approximation schemes supported by our platform. Our choice of approximations is such that they apply both to problems with exogenous and endogenous uncertainty. A decision tree describing the main options available to the user of our platform based on the structure of their problem is provided in Figure~\ref{fig:decision_tree}. Extensions are available in \ROCPP{} which can cater for more classes of problems, see Section~\ref{sec:extensions}.

%\notepv{what to do about multi-stage with constraint uncertainty that cannot be solved using K-adaptability?}

%%%%%%%%%%%%%%%%%%%%%%%%%%%%%%%%%%%%%%%%%%%%%%%%%%%%%%%%%%%%%%%%%%%%%%%%%%%%%%%%%%%%%%
%%%%%%%%%%%%%%%%%%%%%%%%%%%%%%%%%%%%%%%%%%%%%%%%%%%%%%%%%%%%%%%%%%%%%%%%%%%%%%%%%%%%%%

\subsection{Interpretable Decision Rules}
\label{sec:decision_rules}

\paragraph{Constant Decision Rule and Linear Decision Rule.} The most crude (and perhaps most interpretable) decision rules that are available in \ROCPP{} are the \emph{constant decision rule} (CDR) and the \emph{linear decision rule} (LDR), see \cite{BenTal_Book,kuhn_primal_dual_rules}. These apply to binary and real-valued decision variables, respectively. Under the constant decision rule, the binary decisions ${\bm z}_t(\cdot)$ and ${\bm w}_t(\cdot)$ are no longer allowed to adapt to the history of observations -- it is assumed that the decision-maker will take the same action, independently of the realization of the uncertain parameters. Mathematically, we have
$$
{\bm z}_t({\bm \xi}) = {\bm z}_t \text{ and } {\bm w}_t({\bm \xi}) = {\bm w}_t \quad \forall t \in \sets T, \; \forall {\bm \xi} \in \Xi,
$$
for some vectors ${\bm z}_t \in \{0,1\}^{\ell_t}$ and ${\bm w}_t \in \{0,1\}^{k}$, $t \in \sets T$. Under the linear decision rule, the real-valued decisions are modelled as affine functions of the history of observations, i.e.,
$$
{\bm y}_t({\bm \xi}) = {\bm Y}_t {\bm \xi} + {\bm y}_t \quad \forall t \in \sets T,
$$
for some matrix ${\bm Y}_t \in \reals^{n_t \times k}$ and vector ${\bm y}_t \in \reals^{n_t}$. The LDR model leads to very interpretable decisions -- we can think of this decision rule as a scoring rule which assigns different values (coefficients) to each uncertain parameter. These coefficients can be interpreted as the sensitivity of the decision variables to changes in the uncertain parameters. Under the CDR and LDR approximations the adaptive variables in the problem are eliminated and the quantities ${\bm z}_t$, ${\bm w}_t$, ${\bm Y}_t$, and ${\bm y}_t$ become the new decision variables of the problem. 

% \paragraph{Linear Decision Rule.} \emph{Linear decision rules} (LDR) are available in \ROCPP{} to help mitigate the conservatism of CDR for real-valued decision variables. In linear decision rules, the decisions are modelled as linear functions of the history of observations. Mathematically, we have
% $$
% {\bm y}_t({\bm \xi}) = {\bm Y}_t {\bm \xi} + {\bm y}_t \quad \forall t \in \sets T
% $$
% for some matrix ${\bm Y}_t \in \reals^{n_t \times k}$ and vector ${\bm y}_t \in \reals^{n_t}$, that become the new decision variables of the problem. This model leads to very interpretable decisions -- we can think of decision rule as a scoring rule which assigns different values (coefficients) to each uncertain parameter. These coefficients can be interpreted as the sensitivity of the decision variables to changes in the uncertain parameters.% For example, in an inventory management problem, we expect that the amount of an item ordered should increase 
%The sign of these coefficients indicates whether the decision should increase or decrease in response to an increase in the realization of the uncertain parameters and the absolute value of the coefficient (e.g., in an inventory management problem, the amount of an item ordered in response to an increase in observed demand).

\paragraph{Piecewise Constant and Piecewise Linear Decision Rules.} In \emph{piecewise constant} (PWC) and \emph{piecewise linear} (PWL) decision rules, the binary (resp.\ real-valued) adjustable decisions are approximated by functions that are piecewise constant (resp.\ piecewise linear) on a preselected partition of the uncertainty set. Specifically, the uncertainty set $\Xi$ is split into hyper-rectangles of the form 
$
\Xi_{\bm s} \; := \; \left\{ {\bm \xi}\in \Xi \; : \; {\bm c}_{{\bm s}_i-1}^i \leq {\bm \xi}_i < {\bm c}_{{\bm s}_i}^i, \; i=1,\ldots,k \right\},
$ 
where ${\bm s} \in \sets S := \times_{i=1}^{k} \{1,\ldots,{\bm r}_i\} \subseteq \integers^{k}$ and
${\bm c}_1^i \; < \; {\bm c}_2^i \; < \; \cdots \; < \; {\bm c}_{{\bm r}_i-1}^i$, $i=1,\ldots,k$
represent ${\bm r}_i-1$ breakpoints along the ${\bm \xi}_i$ axis. Mathematically, the binary and real-valued decisions are expressible as
$$
{\bm z}_t({\bm \xi}) \; = \; \sum_{{\bm s}\in \sets S} \I{ {\bm \xi} \in  \Xi_{\bm s} }{\bm z}_t^{\bm s}, \quad
{\bm w}_t({\bm \xi}) \; = \; \sum_{{\bm s}\in \sets S} \I{ {\bm \xi} \in  \Xi_{\bm s} }{\bm w}_t^{\bm s},
$$% \quad
%\text{ and } \quad
%
$$
\text{ and } \quad {\bm y}_t({\bm \xi}) \; = \; \sum_{{\bm s}\in \sets S} \I{ {\bm \xi} \in  \Xi_{\bm s} } ( {\bm Y}_t^{\bm s}{\bm \xi} + {\bm y}_t^{\bm s} ), 
$$
for some vectors ${\bm z}_t^{\bm s} \in \reals^{\ell_t}$, ${\bm w}_t^{\bm s} \in \reals^k$, ${\bm y}_t^{\bm s} \in \reals^{n_t}$ and matrices ${\bm Y}_t^{\bm s} \in \reals^{n_t \times k}$, $t\in \sets T$, ${\bm s} \in \sets S$. We can think of this decision rule as a scoring rule which assigns different values (coefficients) to each uncertain parameter; the score assigned to each parameter depends on the subset of the partition in which the uncertain parameter lies. Although less interpretable than CDR/LDR, the PWC/PWL approximation enjoys better optimality properties: it will usually outperform CDR/LDR, since the decisions that can be modelled are more flexible.

The decision rule approximations offered by the \ROCPP{} platform apply to multi-stage problems with both endogenous and exogenous uncertainty. They apply to both stochastic and robust problems, see Figure~\ref{fig:decision_tree}.

%%%%%%%%%%%%%%%%%%%%%%%%%%%%%%%%%%%%%%%%%%%%%%%%%%%%%%%%%%%%%%%%%%%%%%%%%%%%%%%%%%%%%%
%%%%%%%%%%%%%%%%%%%%%%%%%%%%%%%%%%%%%%%%%%%%%%%%%%%%%%%%%%%%%%%%%%%%%%%%%%%%%%%%%%%%%%

\subsection{Contingency Planning via Finite Adaptability}
\label{sec:Kadaptability}

Another solution approach available in \ROCPP{} is the so-called finite adaptability which applies to robust problems with binary decisions, see~\cite{vayanos_ROInfoDiscovery}. It consists in selecting a collection of contingency plans indexed in the set $\sets K := \times_{t\in \sets T} \{ 1 ,\ldots, K_t \}$ today and choosing, at each time $t \in \sets T$, one of the $K_t$ plans for time~$t$ to implement. Mathematically, given $K_t$, $t\in \sets T$, we select ${\bm z}_t^{k_1,\ldots,k_t} \in \{0,1\}^{\ell}_t$ and ${\bm w}_t^{k_1,\ldots,k_t} \in \{0,1\}^{k}$, $k_t \in \{1,\ldots,K_t\}$ for all $t \in \sets T$, in the first stage. At each stage $t \in \sets T$, we select one of the contingency plans $k_t \in \sets K_t$ for implementation. In particular, if up to time $t$, contingency plans $k_1,\ldots,k_t$ have been selected, then ${\bm w}_t^{k_1,\ldots,k_t}$ and ${\bm z}_t^{k_1,\ldots,k_t}$ are implemented at time $t$.

%Another solution approach available in \ROCPP{} is the so-called finite adaptability which applies to robust problems with binary decisions, see~\cite{vayanos_ROInfoDiscovery}. It consists in selecting $K_t \in \naturals$ \emph{candidate} strategies (contingency plans) for each stage $t \in \sets T$ in the first stage and implementing the best of these $K_t$ strategies once the uncertain parameters ${\bm \xi}_i$ such that ${\bm w}_{t-1,i}({\bm \xi}) = 1$, $i = 1,\ldots,k$, are revealed. Mathematically, we select ${\bm z}_t^{k_1,\ldots,k_t} \in \{0,1\}^{\ell}_t$ and ${\bm w}_t^{k_1,\ldots,k_t} \in \{0,1\}^{k}$, $k_t \in \sets K_t:=\{1,\ldots,K_t\}$ for all $t \in \sets T$, in the first stage. At each stage $t \in \sets T$, we select one of the contingency plans $k_t \in \sets K_t$ for implementation. In particular, if up to time $t$, contingency plans $k_1,\ldots,k_t$ have been selected, then ${\bm w}_t^{k_1,\ldots,k_t}$ and ${\bm z}_t^{k_1,\ldots,k_t}$ are implemented at time $t$.

Relative to the piecewise constant decision rule, the finite adaptability approach usually results in better performance in the following sense: the number of contingency plans needed in the finite adaptability approach to achieve a given objective value is never greater than the number of subsets needed in the piecewise constant decision rule to achieve that same value. However, the finite adaptability approximation does not apply to problems with an expectation objective and is thus less flexible in that sense.

%%%%%%%%%%%%%%%%%%%%%%%%%%%%%%%%%%%%%%%%%%%%%%%%%%%%%%%%%%%%%%%%%%%%%%%%%%%%%%%%%%%%%%
%%%%%%%%%%%%%%%%%%%%%%%%%%%%%%%%%%%%%%%%%%%%%%%%%%%%%%%%%%%%%%%%%%%%%%%%%%%%%%%%%%%%%%
%%%%%%%%%%%%%%%%%%%%%%%%%%%%%%%%%%%% EXAMPLES %%%%%%%%%%%%%%%%%%%%%%%%%%%%%%%%%%%%%%%%
%%%%%%%%%%%%%%%%%%%%%%%%%%%%%%%%%%%%%%%%%%%%%%%%%%%%%%%%%%%%%%%%%%%%%%%%%%%%%%%%%%%%%%
%%%%%%%%%%%%%%%%%%%%%%%%%%%%%%%%%%%%%%%%%%%%%%%%%%%%%%%%%%%%%%%%%%%%%%%%%%%%%%%%%%%%%%

\section{Modelling and Solving Decision-Making Problems in \ROCPP}
\label{sec:numerical_results}

In this section, we discuss several robust/stochastic optimization problems, present their mathematical formulations, and discuss how these can be naturally expressed in \ROCPP. \newpv{We also illustrate how optimal solutions can be displayed using our platform.} 

%%%%%%%%%%%%%%%%%%%%%%%%%%%%%%%%%%%%%%%%%%%%%%%%%%%%%%%%%%%%%%%%%%%%%%%%%%%%%%%%%%%%%%
%%%%%%%%%%%%%%%%%%%%%%%%%%%%%%%%%%%%%%%%%%%%%%%%%%%%%%%%%%%%%%%%%%%%%%%%%%%%%%%%%%%%%%
%%%%%%%%%%%%%%%%%%%%%%%%%%%%%%%%%%%%%%%% RSFC %%%%%%%%%%%%%%%%%%%%%%%%%%%%%%%%%%%%%%%%
%%%%%%%%%%%%%%%%%%%%%%%%%%%%%%%%%%%%%%%%%%%%%%%%%%%%%%%%%%%%%%%%%%%%%%%%%%%%%%%%%%%%%%
%%%%%%%%%%%%%%%%%%%%%%%%%%%%%%%%%%%%%%%%%%%%%%%%%%%%%%%%%%%%%%%%%%%%%%%%%%%%%%%%%%%%%%

\subsection{Retailer-Supplier Flexible Commitment Contracts (RSFC)}
\label{sec:retailer}

We model the two-echelon, multiperiod supply chain problem, known as the retailer-supplier flexible commitment (RSFC) problem from~\cite{BenTal_RSFC} in \ROCPP. 

%%%%%%%%%%%%%%%%%%%%%%%%%%%%%%%%%%%%%%%%%%%%%%%%%%%%%%%%%%%%%%%%%%%%%%%%%%%%%%%%%%%%%%
%%%%%%%%%%%%%%%%%%%%%%%%%%%%%%%%%%%%%%%%%%%%%%%%%%%%%%%%%%%%%%%%%%%%%%%%%%%%%%%%%%%%%%

\subsubsection{Problem Description}
\label{sec:RSFC_description}

We consider a finite planning horizon of $T$ periods, $\mathcal T:=\{1,\ldots,T\}$. At the end of each period $t \in \sets T$, the demand ${\bm \xi}_t \in \reals_+$ for the product faced during that period is revealed. We collect the demands for all periods in the vector ${\bm \xi}:=({\bm \xi}_1,\ldots,{\bm \xi}_T)$. We assume that the demand is known to belong to the box uncertainty set
$$
\Xi := \left\{ {\bm \xi} \in \reals^{T} \; : \; {\bm \xi} \in [ \overline{\bm \xi} (1-\rho), \overline {\bm \xi} (1+\rho) ]  \right\},
$$
where $\overline{\bm \xi} := \1 \overline \xi$, $\overline \xi$ is the nominal demand, and $\rho \in [0,1]$ is the budget of uncertainty.

As the time of revelation of the uncertain parameters is exogenous, the information base, encoded in the vectors ${\bm w}_t \in \{0,1\}^T$, $t=0,\ldots,T$, is an input of the problem (data). In particular, it is defined through ${\bm w}_0:={\bm 0}$ and ${\bm w}_{t} := \sum_{\tau=1}^{t} \1_\tau$ for each $t \in \sets T$: the information base for time $t+1$ only contains the demand for the previous periods $\tau \in \{1,\ldots,t\}$.

At the beginning of the planning horizon, the retailer holds an inventory ${\bm x}_1^{\rm i}$ of the product (assumed to be known). At that point, they must specify their commitments ${\bm y}^{\rm{c}} = ( {\bm y}_1^{\rm{c}},\ldots,{\bm y}_{T}^{\rm{c}})$, where ${\bm y}_t^{\rm{c}} \in \reals_+$ represents the amount of the product that they forecast to order at the beginning of time $t \in \sets T$ from the supplier. A penalty cost ${\bm c}^{\rm{dc}+}_t$ (resp.\ ${\bm c}^{\rm dc-}_t$) is incurred for each unit of increase (resp.\ decrease) between the amounts committed for times $t$ and $t-1$. The amount committed for the last order before the beginning of the planning horizon is given by ${\bm y}_0^{\rm c}$. At the beginning of each period, the retailer orders a quantity ${\bm y}_t^{\rm o}({\bm \xi}) \in \reals_+$ from the supplier at unit cost ${\bm c}_t^{\rm{o}}$. These orders may deviate from the commitments made at the beginning of the planning horizon; in this case, a cost ${\bm c}^{\rm dp+}_t$ (resp.\ ${\bm c}^{\rm dp-}_t$) is incurred for each unit that orders ${\bm y}_t^{\rm o}({\bm \xi})$ overshoot (resp.\ undershoot) the plan ${\bm y}_t^{\rm c}$. Given this notation, the inventory of the 
retailer at time $t+1$, $t \in \sets T$, is expressible as
$$
{\bm x}^{\rm i}_{t+1}({\bm \xi}) \; = \; {\bm x}^{\rm i}_t({\bm \xi}) + {\bm y}^{\rm o}_t({\bm \xi}) - {\bm \xi}_{t}.
$$
A holding cost ${\bm c}^{\rm h}_{t+1}$ is incurred for each unit of inventory held on hand at time $t+1$, $t \in \sets T$. A shortage cost ${\bm c}^{\rm s}_{t+1}$ is incurred for each unit of demand lost at time $t+1$, $t \in \sets T$.

The amounts of the product that can be ordered in any given period are constrained by lower and upper bounds denoted by $\underline{\bm y}_t^{\rm o}$ and $\overline {\bm y}_t^{\rm o}$, respectively. Similarly, cumulative orders up to and including time $t \in \sets T$ are constrained to lie in the range $\underline{\bm y}_t^{\rm co}$ to $\overline{\bm y}_t^{\rm co}$. Thus, we have%for each ${\bm \xi} \in \Xi$ the orders for time $t$ are constrained as
$$
\underline{\bm y}_t^{\rm o} \leq {\bm y}_t^{\rm o}({\bm \xi}) \leq \overline{\bm y}_t^{\rm o} \quad  \text{ and } \quad
\underline{\bm y}_t^{\rm co} \leq \sum_{\tau = 1}^t {\bm y}_\tau^{\rm o}({\bm \xi}) \leq \overline{\bm y}_t^{\rm co}.
$$

The objective of the retailer is to minimize their worst-case (maximum) costs. We introduce three sets of auxiliary variables used to linearize the objective function. For each $t\in \sets T$, we let ${\bm y}_t^{\rm dc}$ represent the smallest number that exceeds the costs of deviating from commitments, i.e.,
$$
        {\bm y}_t^{\rm dc} \geq {\bm c}_t^{\rm dc+} ( {\bm y}_t^{\rm c} - {\bm y}_{t-1}^{\rm c} ) %
        \quad \text{ and } \quad
        {\bm y}_t^{\rm dc} \geq {\bm c}_t^{\rm dc-} ({\bm y}_{t-1}^{\rm c} - {\bm y}_{t}^{\rm c} ).
$$
For each $t\in \sets T$ and ${\bm \xi} \in \Xi$, we denote by ${\bm y}_t^{\rm dp}({\bm \xi})$ the smallest number that exceeds the deviations from the plan for time $t$, i.e.,
$$
{\bm y}_t^{\rm dp}({\bm \xi}) \geq {\bm c}_t^{\rm dp+} ({\bm y}_t^{\rm o}({\bm \xi}) - {\bm y}_{t}^{\rm c}) %\\
        \quad \text{ and } \quad
        {\bm y}_t^{\rm dp}({\bm \xi}) \geq {\bm c}_t^{\rm dp-} ({\bm y}_t^{\rm c} - {\bm y}_{t}^{\rm o}({\bm \xi}))
$$
Similarly, for each $t\in \sets T$ and ${\bm \xi} \in \Xi$, we denote by ${\bm y}_{t+1}^{\rm hs}({\bm \xi})$ the smallest number that exceeds the overall holding and shortage costs at time $t+1$, i.e.,
$$
    {\bm y}_{t+1}^{\rm hs}({\bm \xi}) \geq {\bm c}_{t+1}^{\rm h}{\bm x}^{\rm i}_{t+1}({\bm \xi})   %\\
        \quad \text{ and } \quad
    {\bm y}_{t+1}^{\rm hs}({\bm \xi}) \geq -{\bm c}_{t+1}^{\rm s} {\bm x}^{\rm i}_{t+1}({\bm \xi}).
$$
The objective of the retailer is then expressible compactly as
$$
\min \;\; \max_{ {\bm \xi} \in \Xi } \;\; \sum_{t \in \sets T} \; {\bm c}_t^{\rm o} {\bm y}_{t}^{\rm o}({\bm \xi}) +  {\bm y}_t^{\rm dc}({\bm \xi}) + {\bm y}_t^{\rm dp}({\bm \xi}) + {\bm y}_{t+1}^{\rm hs}({\bm \xi}).
$$

The full model for this problem can be found in Electronic Companion~\ref{sec:EC_RSFC_formulation}.

\subsubsection{Model in \ROCPP}

We now present the \ROCPP{} model of the RSFC problem. We assume that the data of the problem have been defined in \cpluspluslogo. The \cpluspluslogo{} variables associated with the problem data are detailed in Table~\ref{tab:RSFCparams}. For example, the lower bounds on the orders $\underline{\bm y}_t^{\rm o}$, $t \in \sets T$, are stored in the map \texttt{OrderLB} which maps each time period to the double representing the lower bound for that period. We discuss how to construct the key elements of the problem here. The full code can be found in Electronic Companion~\ref{sec:EC_RSFC_code}.

\begin{table}[h!]
    \centering\renewcommand{\arraystretch}{\pvtablestretch}
    %{\small{
    \begin{tabular}{c|c|c|c}
    \hline
        Parameter/Index & \cpluspluslogo{} Name  &  \cpluspluslogo{} Type & \cpluspluslogo{} Map Keys \\
        \hline \hline
         $T$ ($t$) & \texttt{T} (\texttt{t}) & \texttt{uint} & NA \\
         %$t$ & \texttt{t} & \texttt{uint} & NA \\
         ${\bm x}_1^{\rm i}$ &  \texttt{InitInventory} & \texttt{double}& NA  \\
         ${\bm y}_0^{\rm c}$ &  \texttt{InitCommit} & \texttt{double} & NA \\
         $\overline \xi$ & \texttt{NomDemand} & \texttt{double}& NA  \\
         $\rho$ & \texttt{rho} & \texttt{double} & NA  \\
         $\underline{\bm y}_t^{\rm o}$, $t\in \sets T$ &  \texttt{OrderLB}  & \texttt{map<uint,double>} & \texttt{t=1\ldots T} \\
         $\overline{\bm y}_t^{\rm o}$, $t\in \sets T$ &  \texttt{OrderUB}  & \texttt{map<uint,double>} & \texttt{t=1\ldots T} \\
         $\underline{\bm y}_t^{\rm co}$, $t\in \sets T$ &  \texttt{CumOrderLB} & \texttt{map<uint,double>}& \texttt{t=1\ldots T}  \\
         $\overline{\bm y}_t^{\rm co}$, $t\in \sets T$ &  \texttt{CumOrderUB} & \texttt{map<uint,double>}& \texttt{t=1\ldots T}  \\
         ${\bm c}^{\rm o}_{t}$, $t\in \sets T$ &  \texttt{OrderCost} & \texttt{map<uint,double>} & \texttt{t=1\ldots T} \\
         ${\bm c}^{\rm h}_{t+1}$, $t\in \sets T$ &  \texttt{HoldingCost} & \texttt{map<uint,double>} & \texttt{t=2\ldots T+1} \\
         ${\bm c}^{\rm s}_{t+1}$, $t\in \sets T$ &  \texttt{ShortageCost} & \texttt{map<uint,double>} & \texttt{t=2\ldots T+1} \\
         ${\bm c}^{\rm dc+}_{t}$, $t\in \sets T$ &  \texttt{CostDCp} & \texttt{map<uint,double>}& \texttt{t=1\ldots T} \\
         ${\bm c}^{\rm dc-}_{t}$, $t\in \sets T$ &  \texttt{CostDCm} & \texttt{map<uint,double>} & \texttt{t=1\ldots T} \\
         ${\bm c}^{\rm dp+}_{t}$, $t\in \sets T$ &  \texttt{CostDPp} & \texttt{map<uint,double>}& \texttt{t=1\ldots T}  \\
         ${\bm c}^{\rm dp-}_{t}$, $t\in \sets T$ &  \texttt{CostDPm} & \texttt{map<uint,double>} & \texttt{t=1\ldots T} \\
        \hline
    \end{tabular}%}}
    \caption{List of model parameters and their associated \cpluspluslogo{} variables for the RSFC problem.}
    \label{tab:RSFCparams}
\end{table}

The RSFC problem is a multi-stage robust optimization problem involving only exogenous uncertain parameters. We begin by creating a model, \texttt{RSFCModel}, in \ROCPP{} that will contain our formulation. All models are pointers to the interface class \texttt{ROCPPOptModelIF}. In this case, we instantiate an object of type \texttt{ROCPPUncOptModel} which is derived from \texttt{ROCPPOptModelIF} and which can model multi-stage optimization problems affected by exogenous uncertain parameters only. The first parameter of the \texttt{ROCPPUncOptModel} constructor corresponds to the maximum time period of any decision variable or uncertain parameter in the problem: in this case, $T+1$. The second parameter of the \texttt{ROCPPUncOptModel} constructor is of the enumerated type \texttt{uncOptModelObjType} which admits two possible values: \texttt{robust}, which indicates a min-max objective; and, \texttt{stochastic}, which indicates an expectation objective. The robust RSFC problem can be initialized as:
\begin{lstlisting}[firstline=1,lastline=2]
// Create an empty robust model with T + 1 periods for the RSFC problem
ROCPPOptModelIF_Ptr RSFCModel(new ROCPPUncOptModel(T+1, robust));
\end{lstlisting}
% \begin{lstlisting}
% // Create an empty robust model with T + 1 periods for the RSFC problem
% ROCPPOptModelIF_Ptr RSFCModel(new ROCPPUncOptModel(T+1, robust));
% \end{lstlisting}
%
We note that in \ROCPP{} all optimization problems are minimization problems.

Next, we discuss how to create the \ROCPP{} variables associated with uncertain parameters and decision variables in the problem. The correspondence between variables is summarized in Table~\ref{tab:RSFCdvs} for convenience.

\begin{table}[ht!]
    \centering \renewcommand{\arraystretch}{\pvtablestretch}
    \begin{tabular}{c|c|c|c}
    \hline
        Variable & \cpluspluslogo{} Name  &  \cpluspluslogo{} Type & \cpluspluslogo{} Map Keys \\
        \hline \hline
         ${\bm \xi}_t$, $t\in \sets T$ &  \texttt{Demand}  & \texttt{map<uint,ROCPPUnc\_Ptr>} & \texttt{t=1\ldots T} \\
         ${\bm y}_t^{\rm o}$, $t\in \sets T$ &  \texttt{Orders}  & \texttt{map<uint,ROCPPVarIF\_Ptr>} & \texttt{t=1\ldots T} \\
         ${\bm y}_t^{\rm c}$, $t\in \sets T$ &  \texttt{Commits}  & \texttt{map<uint,ROCPPVarIF\_Ptr>} & \texttt{t=1\ldots T} \\
         ${\bm y}_t^{\rm dc}$, $t\in \sets T$ &  \texttt{MaxDC}  & \texttt{map<uint,ROCPPVarIF\_Ptr>} & \texttt{t=1\ldots T} \\
         ${\bm y}_t^{\rm dp}$, $t\in \sets T$ &  \texttt{MaxDP}  & \texttt{map<uint,ROCPPVarIF\_Ptr>} & \texttt{t=1\ldots T} \\
         ${\bm y}_{t+1}^{\rm hs}$, $t\in \sets T$ &  \texttt{MaxHS}  & \texttt{map<uint,ROCPPVarIF\_Ptr>} & \texttt{t=2\ldots T+1} \\
         ${\bm x}_{t+1}^{\rm i}$, $t\in \sets T$ &  \texttt{Inventory}  & \texttt{ROCPPExp\_Ptr} & NA \\
         \hline
        \hline
    \end{tabular}
    \caption{List of model variables and uncertainties and their associated \ROCPP{} variables for the RSFC problem.}%{List of decision variables and uncertain parameters in the RSFC problem and their associated \ROCPP{} variables.}
    \label{tab:RSFCdvs}
\end{table}

The uncertain parameters of the RSFC problem are ${\bm \xi}_t$, $t \in \sets T$. We store these in the \texttt{Demand} map, which maps each period to the associated uncertain parameter. Each uncertain parameter is a pointer to an object of type \texttt{ROCPPUnc}. The constructor of the \texttt{ROCPPUnc} class takes two input parameters: the name of the uncertain parameter and the period when that parameter is revealed (first time stage when it is observable).% Thus, the map containing the uncertain parameters can be constructed using the following code:
\lstset{firstnumber=3}
\begin{lstlisting}[firstline=3,lastline=8]
// Create the Demand map to store the uncertain parameters of the problem
map<uint,ROCPPUnc_Ptr> Demand;
// Iterate over all time periods when there is uncertainty
for (uint t=1; t<=T; t++)
    // Create the uncertain parameters, and add them to Demand
    Demand[t] = ROCPPUnc_Ptr(new ROCPPUnc("Demand_"+to_string(t),t+1));
\end{lstlisting}

The main decision variables of the RSFC problem are ${\bm y}_t^c$ and ${\bm y}_t^o$, $t\in \sets T$. The commitment variables ${\bm y}_t^c$ are all static. We store these in the \texttt{Commits} map which maps each time period to the associated commitment decision. The order variables ${\bm y}_t^o$ are allowed to adapt to the history of observed demand realizations. We store these in the \texttt{Orders} map which maps the time period at which the decision is made to the order decision for that period. In \ROCPP, the decision variables are pointers to objects of type \texttt{ROCPPVarIF}. Real-valued static (adaptive) decision variables are modelled using objects of type \texttt{ROCPPStaticVarReal} (\texttt{ROCPPAdaptVarReal}). The constructor of \texttt{ROCPPStaticVarReal} admits three input parameters: the name of the variable, its lower bound, and its upper bound. The constructor of \texttt{ROCPPAdaptVarReal} admits four input parameters: the name of the variable, the time period when the decision is made, and the lower and upper bounds.% The maps containing the order and commitment variables can be constructed using the following code:
\lstset{firstnumber=9}
\begin{lstlisting}[firstline=9,lastline=19]
// Create maps to store the decision variables of the problem
map<uint,ROCPPVarIF_Ptr> Orders, Commits;    // Quantity ordered, Commitments made
// Iterate over all time periods from 1 to T
for (uint t=1; t<=T; t++) {
    // Create the commitment variables (these are static)
    Commits[t]=ROCPPVarIF_Ptr(new ROCPPStaticVarReal("Commit_"+to_string(t),0.));
    if (t==1) // In the first period, the order variables are static
        Orders[t] = ROCPPVarIF_Ptr(new ROCPPStaticVarReal("Order_"+to_string(t),OrderLB[t],OrderUB[t]));
    else // In the other periods, the order variables are adaptive
        Orders[t] = ROCPPVarIF_Ptr(new ROCPPAdaptVarReal("Order_"+to_string(t),t,OrderLB[t],OrderUB[t]));
}
\end{lstlisting}

The RSFC problem also involves the auxiliary variables ${\bm y}_t^{\rm dc}$, ${\bm y}_t^{\rm dp}$, and ${\bm y}_{t+1}^{\rm hs}$, $t\in \sets T$. We store the ${\bm y}_t^{\rm dc}$ variables in the map \texttt{MaxDC}. These variables are all static. We store the ${\bm y}_t^{\rm dp}$ variables in the map \texttt{MaxDP}. We store the ${\bm y}_t^{\rm hs}$ variables in the map  \texttt{MaxHS}. Since the orders placed and inventories held change based on the demand realization, the variables stored in \texttt{MaxDP} and \texttt{MaxHS} are allowed to adapt to the history of observations. All maps map the index of the decision to the actual decision variable. The procedure to build these maps exactly parallels the approach above and we thus omit it. We refer the reader to lines \ref{ln:retailer_vars_2_start}-\ref{ln:retailer_vars_2_end} in Section~\ref{sec:EC_RSFC_code} for the code to build these maps. % The maps containing these auxiliary variables can be constructed using the following code:
% %
% \lstset{firstnumber=20}
% \lstinputlisting[firstline=20,lastline=34]{retailer.cc}

Having defined the model, the uncertain parameters, and the decision variables of the problem, we are now ready to formulate the constraints. To express our constraints in an intuitive fashion, we create an expression variable (type \texttt{ROCPPExpr}), which we call \texttt{Inventory} that stores the amount of inventory held at the beginning of each period. This is computed by adding to the initial inventory \texttt{InitInventory} the amount ordered at each period and subtracting the demand faced. Similarly, we create an \texttt{ROCPPExpr} to store the cumulative orders placed. This is obtained by adding orders placed at each period. Constraints can be created using the operators ``\lstinline{<=}'', ``\lstinline{>=}'', or ``\lstinline{==}'' and added to the problem using the \texttt{add\_constraint()} function. We show how to construct the cumulative order constraints and the lower bounds on the shortage and holding costs. The code to add the remaining constraints can be found in lines~\ref{ln:retailer_extra_cstr_start}-\ref{ln:retailer_extra_cstr_end} of Section~\ref{sec:EC_RSFC_code}.
\lstset{firstnumber=35}
\begin{lstlisting}[firstline=35,lastline=52]
// Create the constraints of the problem
// Create an expression for the amount of inventory held and initialize it
ROCPPExpr_Ptr Inventory(new ROCPPExpr());
Inventory = Inventory + InitInventory;
// Create an expression for the cumulative amount ordered
ROCPPExpr_Ptr CumOrders(new ROCPPExpr());
// Iterate over all time periods and add the constraints to the problem
for (uint t=1; t<=T; t++) {
    // Create the upper and lower bounds on the cumulative orders
    CumOrders = CumOrders + Orders[t];
    RSFCModel->add_constraint(CumOrders >= CumOrderLB[t]);
    RSFCModel->add_constraint(CumOrders <= CumOrderUB[t]);
    // Update the inventory
    Inventory = Inventory + Orders[t] - Demand[t];
    // Create upper bound on shortage/holding costs
    RSFCModel->add_constraint(MaxHS[t+1] >= HoldingCost[t+1]*Inventory);
    RSFCModel->add_constraint(MaxHS[t+1] >= (-1.*ShortageCost[t+1]*Inventory));
}
\end{lstlisting}

The objective function of the RSFC problem consists in minimizing the sum of all costs over time. We create the \texttt{ROCPPExpr} expression \texttt{RSFCObj} to which we add all terms by iterating over time. We then set \texttt{RSFCObj} as the objective function of the \texttt{RSFCModel} model by using the \texttt{set\_objective()} function.
\lstset{firstnumber=68}
\begin{lstlisting}[firstline=68,lastline=73]
// Create an expression that will contain the objective function
ROCPPExpr_Ptr RSFCObj(new ROCPPExpr());
// Iterate over all periods and add the terms to the objective function
for (uint t=1; t<=T; t++)
    RSFCObj = RSFCObj + OrderCost*Orders[t] + MaxDC[t] + MaxDP[t] + MaxHS[t+1];
RSFCModel->set_objective(RSFCObj); // Add the objective to the problem
\end{lstlisting}

Finally, we create a box uncertainty set for the demand.
\lstset{firstnumber=74}
\begin{lstlisting}[firstline=74,lastline=78]
for (uint t=1; t<=T; t++) |\label{ln:uset_start}|{
    // Add the upper and lower bounds on the demand to the uncertainty set
    RSFCModel->add_constraint_uncset(Demand[t] >= NomDemand*(1.0-rho));
    RSFCModel->add_constraint_uncset(Demand[t] <= NomDemand*(1.0+rho));
} |\label{ln:uset_end}|
\end{lstlisting}

Having formulated the RSFC problem in \ROCPP, we turn to solving it.

%%%%%%%%%%%%%%%%%%%%%%%%%%%%%%%%%%%%%%%%%%%%%%%%%%%%%%%%%%%%%%%%%%%%%%%%%%%%%%%%%%%%%%
%%%%%%%%%%%%%%%%%%%%%%%%%%%%%%%%%%%%%%%%%%%%%%%%%%%%%%%%%%%%%%%%%%%%%%%%%%%%%%%%%%%%%%

\subsubsection{Solution in \ROCPP}
\label{sec:RSFC_solution}

From~\cite{BenTal_RSFC}, LDRs are optimal in this case. Thus, it suffices to approximate the real-valued adaptive variables in the problem by linear decision rules, then robustify the problem using duality theory, and finally solve it using an off-the-shelf deterministic solver. This process is streamlined in \ROCPP{}.
\lstset{firstnumber=79}
\begin{lstlisting}[firstline=79,lastline=86]
// Construct the linear/constant decision rule approximator
ROCPPApproximator_Ptr pLDRApprox(new ROCPPLCDRApprox(RSFCModel));
// Approximate the adaptive decisions using the LDR/CDR approximator 
ROCPPMISOCP_Ptr RSFCModelLDR(pLDRApprox->DoMyThing(RSFCModel));
// Construct the solver (in this case, use gurobi as deterministic solver)
ROCPPSolver_Ptr pSolver(new ROCPPGurobi(SolverParams()));
// Solve the problem
pSolver->solve(RSFCModelLDR);
\end{lstlisting}
We consider the instance of RSFC detailed in Electronic Companion~\ref{sec:EC_RSFC_instance}. The following output is displayed when executing the above code on this instance.
\begin{lstlisting}[numbers=none,language=bash]
=========================================================================
=================== APPROXIMATING USING LDR AND CDR ===================== 
========================================================================= 
11 of 119 constraints robustified
...
110 of 119 constraints robustified
Total time to approximate and robustify: 0 seconds
=========================================================================
\end{lstlisting}
This states that the problem has 119 constraints in total and that the time it took to approximate it and obtain its robust counterpart was under half a second. % After this output is shown, the output from the solver is displayed. We omit this here.
Next, we showcase how optimal solutions to the problem can be retrieved in \ROCPP.
\lstset{firstnumber=87}
\begin{lstlisting}[firstline=87,lastline=90]
// Retrieve the optimal solution from the solver
map<string,double> optimalSln(pSolver->getSolution());
// Print the optimal decision (from the original model)
pLDRApprox->printOut(RSFCModelLDR, optimalSln, Orders[10]));
\end{lstlisting}
The following output is displayed when executing the above code.
\begin{lstlisting}[numbers=none,language=bash]
Order_10 =  + 0*Demand_1 + 0*Demand_2 + 0*Demand_3 + 0*Demand_4 + 0*Demand_5 + 0*Demand_6 + 0*Demand_7 + 0*Demand_8 + 1*Demand_9 - 0.794
\end{lstlisting}
Thus, the optimal linear decision rule for the amount to order at stage 10 is ${\bm y}^{\rm o}_{10} ({\bm \xi}) = {\bm \xi}_9 - 0.794$ for this specific instance.
To get the optimal objective value of \texttt{RSFCModelLDR}, we can use the following command, which returns $13531.7$ in this instance.
\lstset{firstnumber=92}
\begin{lstlisting}[firstline=92,lastline=92]
double optVal(pSolver->getOptValue());
\end{lstlisting}

%%%%%%%%%%%%%%%%%%%%%%%%%%%%%%%%%%%%%%%%%%%%%%%%%%%%%%%%%%%%%%%%%%%%%%%%%%%%%%%%%%%%%%
%%%%%%%%%%%%%%%%%%%%%%%%%%%%%%%%%%%%%%%%%%%%%%%%%%%%%%%%%%%%%%%%%%%%%%%%%%%%%%%%%%%%%%

\subsubsection{Variant: Ellipsoidal Uncertainty Set}

In~\cite{BenTal_RSFC}, the authors also investigated ellipsoidal uncertainty sets for the demand. These take the form
$$
\Xi := \{ {\bm \xi} \in \reals_+^{T} : \| {\bm \xi} - \overline {\bm \xi} \|_2 \leq \Omega  \},
$$
where $\Omega$ is a safety parameter. Letting \texttt{Omega} represent $\Omega$, this ellipsoidal uncertainty set can be used instead of the box uncertainty set by replacing lines \ref{ln:uset_start}-\ref{ln:uset_end} in the \ROCPP{} code for the RSFC problem with the following code:
\begin{lstlisting}[numbers=none]
// Create a vector that will contain all the elements of the norm term
vector<ROCPPExpr_Ptr> EllipsoidElements;
// Populate the vector with the difference between demand and nominal demand
for (uint t=1; t<=T; t++)
    EllipsoidElements.push_back(Demand[t+1] - NominalDemand);
// Create the norm term
boost::shared_ptr<ConstraintTermIF> EllipsTerm(new NormTerm(EllipsoidElements));
// Create the ellipsoidal uncertainty constraint
RSFCModel->add_constraint_uncset(EllipsTerm <= Omega);
\end{lstlisting}
The solution approach from Section~\ref{sec:RSFC_solution} applies as is with this ellipsoidal set. The time it takes to robustify the problem is again under half a second. In this case, the optimal objective value under LDRs is $14,814.3$. The optimal linear decision rule is given by:
\begin{lstlisting}[numbers=none,language=bash]
Order_10 =  + 0.0305728*Demand_1 + 0.0567*Demand_2 + 0.0739*Demand_3 + 0.0887*Demand_4 + 0.101*Demand_5 + 0.115*Demand_6 + 0.142*Demand_7 + 0.179*Demand_8 + 0.231*Demand_9 - 3.33
\end{lstlisting}
%

%%%%%%%%%%%%%%%%%%%%%%%%%%%%%%%%%%%%%%%%%%%%%%%%%%%%%%%%%%%%%%%%%%%%%%%%%%%%%%%%%%%%%%
%%%%%%%%%%%%%%%%%%%%%%%%%%%%%%%%%%%%%%%%%%%%%%%%%%%%%%%%%%%%%%%%%%%%%%%%%%%%%%%%%%%%%%
%%%%%%%%%%%%%%%%%%%%%%%%%%%%%%%%%%%%%%%% PB %%%%%%%%%%%%%%%%%%%%%%%%%%%%%%%%%%%%%%%%%%
%%%%%%%%%%%%%%%%%%%%%%%%%%%%%%%%%%%%%%%%%%%%%%%%%%%%%%%%%%%%%%%%%%%%%%%%%%%%%%%%%%%%%%
%%%%%%%%%%%%%%%%%%%%%%%%%%%%%%%%%%%%%%%%%%%%%%%%%%%%%%%%%%%%%%%%%%%%%%%%%%%%%%%%%%%%%%

\subsection{Robust Pandora's Box Problem}
\label{sec:pandora}

We consider a robust variant of the celebrated stochastic Pandora Box (PB) problem due to~\cite{Weitzman_1979}. This problem models selection from a set of unknown, alternative options, when evaluation is costly.% For example, the problem of hiring a skilled worker, where only one hire can be made and where the evaluation of each candidate is an expensive procedure, can be cast as a Pandora's box problem.

%%%%%%%%%%%%%%%%%%%%%%%%%%%%%%%%%%%%%%%%%%%%%%%%%%%%%%%%%%%%%%%%%%%%%%%%%%%%%%%%%%%%%%
%%%%%%%%%%%%%%%%%%%%%%%%%%%%%%%%%%%%%%%%%%%%%%%%%%%%%%%%%%%%%%%%%%%%%%%%%%%%%%%%%%%%%%

\subsubsection{Problem Description}
\label{sec:RPBdescription}

There are~$I$ boxes indexed in~$\mathcal{I}:=\{1,\ldots,I\}$ that we can choose or not to open over the planning horizon $\mathcal T:=\{1,\ldots,T\}$.  Opening box $i \in \sets I$ incurs a cost ${\bm c}_i \in \reals_+$. Each box has an unknown value ${\bm \xi}_i \in \reals$, $i\in \sets I$, which will only be revealed if the box is opened. At the beginning of each time $t \in \sets T$, we can either select a box to open or keep one of the opened boxes, earn its value (discounted by $\theta^{t-1}$), and stop the search.

We assume that the box values are restricted to lie in the set
$$
\Xi := \left\{
\bm{\xi} \in \reals^{I} \; : \; \exists {\bm \zeta} \in [-1, 1]^M , \;
{\bm \xi}_i \; = \; (1+{\bm \Phi}_i^\top {\bm \zeta}/2) \overline{\bm \xi}_i \quad \forall i \in \sets I \right\},
$$
where ${\bm \zeta} \in \reals^M$ represent $M$ risk factors, and $\overline{\bm \xi} \in \reals^I$ collects the nominal box values.

In this problem, the box values are endogenous uncertain parameters whose time of revelation can be controlled by the box open decisions. Thus, the information base, encoded by the vector ${\bm w}_{t}({\bm \xi}) \in \{0,1\}^I$, $t\in \sets T$, is a decision variable. In particular, ${\bm w}_{t,i}({\bm \xi}) = 1$ if and only if box $i \in \sets I$ has been opened on or before time $t \in \sets T$ in scenario ${\bm \xi}$. We assume that ${\bm w}_0({\bm \xi})={\bm 0}$ so that no box is opened before the beginning of the planning horizon. We denote by ${\bm z}_{t,i}({\bm \xi}) \in \{0,1\}$ the decision to keep box $i \in \sets I$ and stop the search at time $t \in \sets T$.

The requirement that at most one box be opened at each time $t\in \sets T$ and that no box be opened if we have stopped the search can be captured by the constraint
\begin{equation}
\sum_{i \in \sets I} ( {\bm w}_{t,i}({\bm \xi}) - {\bm w}_{t-1,i}({\bm \xi}) ) \; \leq \; 1 - \sum_{\tau=1}^t \sum_{i \in \sets I} {\bm z}_{t,i}({\bm \xi}) \quad \forall t \in \sets T.
\label{eq:PB_constraint_1box}
\end{equation}
The requirement that only one of the opened boxes can be kept is expressible as
\begin{equation}
{\bm z}_{t,i}({\bm \xi}) \; \leq \; {\bm w}_{t-1,i}({\bm \xi}) \quad \forall t \in \sets T, \; \forall i \in \sets I.
\label{eq:PB_constraint_keep}
\end{equation}

The objective of the PB problem is to select the sequence of boxes to open and the box to keep so as to maximize worst-case net profit. Since the decision to open box $i$ at time $t$ can be expressed as the difference $({\bm w}_{t,i}-{\bm w}_{t-1,i})$, the objective of the PB problem is
$$
\max \;\; \min_{{\bm \xi} \in \Xi} \;\; \sum_{t \in \sets T} \sum_{i \in \sets I} \theta^{t-1} {\bm \xi}_i {\bm z}_{t,i}({\bm \xi}) - {\bm c}_i ({\bm w}_{t,i}({\bm \xi})-{\bm w}_{t-1,i}({\bm \xi})).
$$
The mathematical model for this problem can be found in Electronic Companion~\ref{sec:EC_PB_formulation}.

%%%%%%%%%%%%%%%%%%%%%%%%%%%%%%%%%%%%%%%%%%%%%%%%%%%%%%%%%%%%%%%%%%%%%%%%%%%%%%%%%%%%%%

\subsubsection{Model in RO\cpluspluslogo}
\label{sec:RPBcode}

We present the \ROCPP{} model for the PB problem. We assume that the data of the problem have been defined in \cpluspluslogo{} as summarized in Table~\ref{tab:PBparams}. 
\begin{table}[ht!]
    \centering \renewcommand{\arraystretch}{0.75}
    \begin{tabular}{c|c|c|c}
    \hline
        Model Parameter & \cpluspluslogo{} Name  &  \cpluspluslogo{} Variable Type & \cpluspluslogo{} Map Keys \\
        \hline \hline
         $\theta$ & \texttt{theta} & \texttt{double} & NA\\
         $T$ ($t$) & \texttt{T} (\texttt{t}) & \texttt{uint} & NA \\
         $I$ ($i$) & \texttt{I} (\texttt{i}) & \texttt{uint} & NA \\
         $M$ ($m$) & \texttt{M} (\texttt{m}) & \texttt{uint} & NA \\
         %$T$ & \texttt{T} & \texttt{uint} & NA \\
         %$t$ & \texttt{t} & \texttt{uint} & NA \\
         %$I$ & \texttt{I} & \texttt{uint} & NA \\
         %$i$ & \texttt{i} & \texttt{uint} & NA \\
         %$M$ & \texttt{M} & \texttt{uint} & NA \\
         %$m$ & \texttt{m} & \texttt{uint} & NA \\
         ${\bm c}_i$,\; $i\in\sets I$ &  \texttt{CostOpen} & \texttt{map<uint,double>}& \texttt{i=1\ldots I}  \\
         $\overline{\bm \xi}_i$,\; $i\in\sets I$ & \texttt{NomVal} & \texttt{map<uint,double>} & \texttt{i=1\ldots I}\\
         ${\bm \Phi}_{im}$, $i\in\sets I$,\;$m\in\sets M$ & \texttt{FactorCoeff} & \texttt{map<uint,map<uint,double> >} & \texttt{i=1\ldots I, m=1\ldots M}\\
        \hline
    \end{tabular}
    \caption{List of model parameters and their associated \cpluspluslogo{} variables for the PB problem.}
    \label{tab:PBparams}
\end{table}
\vspace{-0.5cm}

The PB problem is a multi-stage robust optimization problem involving uncertain parameters whose time of revelation is decision-dependent. Such models can be stored in the \texttt{ROCPPDDUOptModel} class which is derived from \texttt{ROCPPOptModelIF}.% The robust PB problem can be initialized as:
\lstset{firstnumber=1}
\begin{lstlisting}[firstline=1,lastline=2]
// Create an empty robust model with T periods for the PB problem
ROCPPOptModelIF_Ptr PBModel(new ROCPPDDUOptModel(T, robust));
\end{lstlisting}

Next, we create the \ROCPP{} variables associated with uncertain parameters and decision variables in the problem. The correspondence between variables is summarized in Table~\ref{tab:PBvars}.
\vspace{-0.5cm}
\begin{table}[ht!]
    \centering \renewcommand{\arraystretch}{0.75}
    \begin{tabular}{c|c|c|c}
    \hline
        Parameter & \cpluspluslogo{} Nm.  &  \cpluspluslogo{} Type & \cpluspluslogo{} Map Keys \\
        \hline \hline
         ${\bm z}_{t,i}$,\; $i\in\sets I$,\;$t\in\sets T$ &  \texttt{Keep} & \texttt{map<uint,map<uint,ROCPPVarIF\_Ptr> >}& \texttt{1\ldots T, 1\ldots I}  \\
         ${\bm w}_{t,i}$,\; $i\in\sets I$,\;$t\in\sets T$ &  \texttt{MeasVar} &
         \texttt{map<uint,map<uint,ROCPPVarIF\_Ptr> >}& \texttt{1\ldots T, 1\ldots I}  \\
         ${\bm \zeta}_m$,\; $m\in\sets M$ &  \texttt{Factor} &
         \texttt{map<uint,ROCPPUnc\_Ptr>}& \texttt{m=1\ldots M}  \\
         ${\bm \xi}_i$,\; $i\in\sets I$ &  \texttt{Value} &
         \texttt{map<uint,ROCPPUnc\_Ptr>}& \texttt{i=1\ldots I}  \\
        \hline
    \end{tabular}
    \caption{List of model variables and uncertainties and their associated \cpluspluslogo{} variables for the PB problem.}
    \label{tab:PBvars}
\end{table}

The uncertain parameters of the PB problem are ${\bm \xi} \in \reals^I$ and ${\bm \zeta} \in \reals^M$. We store the \ROCPP{} variables associated with these in the \texttt{Value} and \texttt{Factor} maps, respectively. Recall that the first and second (optional) parameters of the \texttt{ROCPPUnc} constructor are the name of the parameter and the time period when it is observed. As~${\bm \xi}$ has a time of revelation that is decision-dependent, we can omit the second parameter when we construct the associated \ROCPP{} variables. The \texttt{ROCPPUnc} constructor also admits a third (optional) parameter with default value \texttt{true} that indicates if the uncertain parameter is observable. As ${\bm \zeta}$ is an auxiliary uncertain parameter, we set its time period as being, e.g., 1 and indicate through the third parameter in the constructor of \texttt{ROCPPUnc} that this parameter is not observable. 
\lstset{firstnumber=3}
\begin{lstlisting}[firstline=3,lastline=10]
// Create empty maps to store the uncertain parameters
map<uint, ROCPPUnc_Ptr> Value, Factor;
for (uint i = 1; i <= I; i++)
    // Create the uncertainty associated with box i and add it to Value
    Value[i] = ROCPPUnc_Ptr(new ROCPPUnc("Value_"+to_string(i)));
for (uint m = 1; m <= M; m++)
    // The risk factors are not observable
    Factor[m]= ROCPPUnc_Ptr(new ROCPPUnc("Factor_"+to_string(m),1,false));
\end{lstlisting}

The decision variables of the problem are the measurement variables~${\bm w}$ and the variables~${\bm z}$ which decide on the box to keep. We store these in the maps \texttt{MeasVar} and \texttt{Keep}, respectively. In \ROCPP, the measurement variables are created automatically for all time periods in the problem by calling the \texttt{add\_ddu()} function which is a public member of \texttt{ROCPPOptModelIF}. This function admits four input parameters: an uncertain parameter, the first and last time period when the decision-maker can choose to observe that parameter, and the cost for observing the parameter. In this problem, cost for observing ${\bm \xi}_i$ is equal to~${\bm c}_i$. The measurement variables constructed in this way can be recovered using the \texttt{getMeasVar()} function which admits as inputs the name of an uncertain parameter and the time period for which we want to recover the measurement variable associated with that uncertain parameter. The boolean \texttt{Keep} variables can be built in \ROCPP{} using the constructors of the \texttt{ROCPPStaticVarBool} and \texttt{ROCPPAdaptVarBool} classes for the static and adaptive variables, respectively, see lines \ref{ln:keep_vars_start}-\ref{ln:keep_vars_end} in Section~\ref{sec:EC_PB_code}. We omit this here.
\lstset{firstnumber=11}
\begin{lstlisting}[firstline=11,lastline=18]
map<uint, map<uint, ROCPPVarIF_Ptr> > MeasVar;
for (uint i = 1; i <= I; i++) {
    // Create the measurement variables associated with the value of box i
    PBModel->add_ddu(Value[i], 1, T, obsCost[i]);
    // Get the measurement variables and store them in MeasVar
    for (uint t = 1; t <= T; t++)
        MeasVar[t][i] = PBModel->getMeasVar(Value[i]->getName(), t);
}
\end{lstlisting}

Having created the decision variables and uncertain parameters, we turn to adding the constraints to the model. To this end, we use the \texttt{StoppedSearch} expression, which tracks the running sum of the \texttt{Keep} variables, to indicate if at any given point in time, we have already decided to keep one box and stop the search. We also use the \texttt{NumOpened} expression which, at each period, stores the expression for the total number of boxes that we choose to open in that period. Since the construction of the constraints is similar to that used in the retailer problem, we omit it here and refer the reader to lines \ref{ln:pandora_cstr_start}-\ref{ln:pandora_cstr_end} in Section~\ref{sec:EC_PB_code}. %Using these expressions, the constraints can be added to the problem using the following code.
% %
% \lstset{firstnumber=27}
% \lstinputlisting[firstline=27,lastline=44]{pandora.cc}

Next, we create the uncertainty set and the objective function.% We note that the uncertainty set accepts equality constraints directly.% which can be done using the following code.% The first set of constraints adds the upper and lower bounds on ${\bm \xi}$, while the second set of constrains adds the linear constraints that express the box values in terms of the risk factors.
\lstset{firstnumber=46}
\begin{lstlisting}[firstline=46,lastline=65]
// Create the uncertainty set constraints and add them to the problem
// Add the upper and lower bounds on the risk factors
for (uint m = 1; m <= M; m++) {
    PBModel->add_constraint_uncset(Factor[m] >= -1.0);
    PBModel->add_constraint_uncset(Factor[m] <= 1.0);
}
// Add the expressions for the box values in terms of the risk factors
for (uint i = 1; i <= I; i++) {
    ROCPPExpr_Ptr ValueExpr(new ROCPPExpr());
    for (uint m = 1; m <= M; m++)
        ValueExpr = ValueExpr + RiskCoeff[i][m]*Factor[m];
    PBModel->add_constraint_uncset(  Value[i] == (1.+0.5*ValueExpr) * NomVal[i] );
}
// Create the objective function expression
ROCPPExpr_Ptr PBObj(new ROCPPExpr()); 
for (uint t = 1; t <= T; t++)
    for (uint i = 1; i <= I; i++)
        PBObj = PBObj + pow(theta,t-1)*Value[i]*Keep[t][i];
// Set objective (multiply by -1 for maximization)
PBModel->set_objective(-1.0*PBObj);
\end{lstlisting}
We emphasize that the observation costs were automatically added to the objective function when we called the \texttt{add\_ddu()} function.

% Finally, we conclude the construction of the PB model by adding the objective function.%, which can be achieved using the following code:
% %
% \lstset{firstnumber=last}
% \begin{lstlisting}

%%%%%%%%%%%%%%%%%%%%%%%%%%%%%%%%%%%%%%%%%%%%%%%%%%%%%%%%%%%%%%%%%%%%%%%%%%%%%%%%%%%%%%

\subsubsection{Solution in \ROCPP}
\label{sec:pandora_solution}

The PB problem is a multi-stage robust problem with decision-dependent information discovery, see~\cite{DDI_VKB, vayanos_ROInfoDiscovery}. \ROCPP{} offers two options for solving this class of problems: finite adaptability and piecewise constant decision rules, see Section~\ref{sec:approximation_schemes}. Here, we illustrate how to solve PB using the finite adaptability approach, see Section~\ref{sec:Kadaptability}. Letting \texttt{uint} \texttt{K} store the number of contingency plans $K$ per period, the process of computing the optimal contingency plans is streamlined in \ROCPP{}.
\lstset{firstnumber=66}
\begin{lstlisting}[firstline=66,lastline=69]
// Construct the finite adaptability approximator
ROCPPApproximator_Ptr pKAdaptaApprox(new ROCPPKadaptApprox(PBModel,K));
// Approximate the decisions using finite adaptability and robustify
ROCPPMISOCP_Ptr PBModelKAdapt(pKAdaptaApprox->DoMyThing(PBModel));
\end{lstlisting}
We consider the instance of PB detailed in Electronic Companion~\ref{sec:EC_PB_instance} for which $T=4$, $M=4$, and $I=5$. For $K=1$ (resp.\ $K=2$ and $K=3$), the problem takes under half a second (resp.\ under half a second and 6 seconds) to approximate and robustify. Its objective value is 2.12 (resp.\ 9.67 and 9.67). Note that with $T=4$ and $K=2$ (resp.\ $K=3$), the total number of contingency plans is $K^{T-1} = 8$ (resp.\ 27).

Next, we showcase how optimal contingency plans can be retrieved in \ROCPP.
\lstset{firstnumber=73}
\begin{lstlisting}[firstline=73,lastline=76]
// Retrieve the optimal solution from the solver
map<string,double> optimalSln(pSolver->getSolution());
// Print the optimal decision (from the original model)
pKAdaptaApprox->printOut(PBModelKAdapt, optimalSln, Keep[4][2]);
\end{lstlisting}
When executing this code, the values of all variables ${\bm w}^{t,k_1 \ldots k_t}$ for all contingency plans $(k_1,\ldots,k_t) \in \times_{\tau=1}^t \sets K^\tau$ are printed. We show here the subset of the output associated with contingency plans where ${\bm z}_{2,4}({\bm \xi})$ equals 1 (for the case $K=2$).
\begin{lstlisting}[numbers=none,language=bash]
Value of variable Keep_4_2 under contingency plan (1-1-2-2) is: 1
\end{lstlisting}
Thus, at time 4, we will keep the second box if and only if the contingency plan we choose is $(k_1,k_2,k_3,k_4)=(1,1,2,2)$. We can display the first time that an uncertain parameter is observed using the following \ROCPP{} code.
\lstset{firstnumber=76}
\begin{lstlisting}[firstline=76,lastline=76]
pKAdaptaApprox->printOut(PBModelKAdapt, optimalSln, Keep[4][2]);
\end{lstlisting}
When executing this code, the time when ${\bm \xi}_4$ is observed under each contingency plan $(k_1,\ldots,k_T) \in \times_{\tau \in \sets T} \sets K^t$ is printed. In this case, part of the output we get is as follows.%If ${\bm \xi}_4$ is not observed at any time under some contingency plan, this is also explicitly stated. 
\begin{lstlisting}[numbers=none,language=bash]
Parameter Value_4 under contingency plan (1-1-2-1) is observed at time 3
Parameter Value_4 under contingency plan (1-2-1-1) is never observed
\end{lstlisting}
Thus, in an optimal solution, ${\bm \xi}_4$ is opened at time 3 under contingency plan $(k_1,k_2,k_3,k_4)=(1,1,2,1)$. On the other hand it is never opened under contingency plan $(1,2,1,1)$.

%%%%%%%%%%%%%%%%%%%%%%%%%%%%%%%%%%%%%%%%%%%%%%%%%%%%%%%%%%%%%%%%%%%%%%%%%%%%%%%%%%%%%%
%%%%%%%%%%%%%%%%%%%%%%%%%%%%%%%%%%%%%%%%%%%%%%%%%%%%%%%%%%%%%%%%%%%%%%%%%%%%%%%%%%%%%%
%%%%%%%%%%%%%%%%%%%%%%%%%%%%%%%%%%%%%%%% BB %%%%%%%%%%%%%%%%%%%%%%%%%%%%%%%%%%%%%%%%%%
%%%%%%%%%%%%%%%%%%%%%%%%%%%%%%%%%%%%%%%%%%%%%%%%%%%%%%%%%%%%%%%%%%%%%%%%%%%%%%%%%%%%%%
%%%%%%%%%%%%%%%%%%%%%%%%%%%%%%%%%%%%%%%%%%%%%%%%%%%%%%%%%%%%%%%%%%%%%%%%%%%%%%%%%%%%%%

\subsection{Stochastic Best Box Problem with Uncertain Observation Costs}
\label{sec:best_box}

%%%%%%%%%%%%%%%%%%%%%%%%%%%%%%%%%%%%%%%%%%%%%%%%%%%%%%%%%%%%%%%%%%%%%%%%%%%%%%%%%%%%%%
%%%%%%%%%%%%%%%%%%%%%%%%%%%%%%%%%%%%%%%%%%%%%%%%%%%%%%%%%%%%%%%%%%%%%%%%%%%%%%%%%%%%%%

We consider a variant of Pandora's Box problem, known as Best Box (BB), in which observation costs are uncertain and subject to a budget constraint. We assume that the decision-maker is interested in maximizing the expected value of the box kept.

\subsubsection{Problem description}
\label{sec:best_box_description}

There are~$I$ boxes indexed in~$\mathcal{I}:=\{1,\ldots,I\}$ that we can choose or not to open over the planning horizon $\mathcal T:=\{1,\ldots,T\}$.  Opening box $i \in \sets I$ incurs an uncertain cost ${\bm \xi}_i^{\rm c} \in \reals_+$. Each box has an unknown value ${\bm \xi}_i^{\rm v} \in \reals$, $i\in \sets I$. The value of each box and the cost of opening it will only be revealed if the box is opened. The total cost of opening boxes cannot exceed budget $B$. At each period $t \in \sets T$, we can either open a box or keep one of the opened boxes, earn its value (discounted by $\theta^{t-1}$), and stop the search.

We assume that box values and costs are uniformly distributed in the set 
$
\Xi := \left\{
{\bm \xi}^{\rm v} \in \reals_+^{I}, \; {\bm \xi}^{\rm c} \in \reals_+^{I} \; : \; {\bm \xi}^{\rm v} \; \leq \; \overline{\bm \xi}^{\rm v}, \; {\bm \xi}^{\rm c} \; \leq \; \overline{\bm \xi}^{\rm c}  \right\}, 
$ where $\overline{\bm \xi}^{\rm v}, \; \overline{\bm \xi}^{\rm c} \in \reals^I$. 

In this problem, the box values and costs are endogenous uncertain parameters whose time of revelation can be controlled by the box open decisions. For each $i\in \sets I$, and $t\in \sets T$, we let ${\bm w}_{t,i}^{\rm v}({\bm \xi})$ and ${\bm w}_{t,i}^{\rm c}({\bm \xi})$ indicate if parameters ${\bm \xi}_i^{\rm v}$ and  ${\bm \xi}_i^{\rm c}$ have been observed on or before time $t$. In particular ${\bm w}_{t,i}^{\rm v}({\bm \xi}) = {\bm w}_{t,i}^{\rm c}({\bm \xi})$ for all $i$, $t$, and ${\bm \xi}$. We assume that ${\bm w}_0({\bm \xi})={\bm 0}$ so that no box is opened before the beginning of the planning horizon. We denote by ${\bm z}_{t,i}({\bm \xi}) \in \{0,1\}$ the decision to keep box $i \in \sets I$ and stop the search at time $t \in \sets T$. The requirement that at most one box be opened at each time $t\in \sets T$ and that no box be opened if we have stopped the search can be captured in a manner that parallels constraint~\eqref{eq:PB_constraint_1box}. Similarly, the requirement that only one of the opened boxes can be kept can be modelled using a constraint similar to~\eqref{eq:PB_constraint_keep}. The budget constraint and objective can be expressed compactly as
$$
\sum_{i \in \sets I} {\bm \xi}_i^{\rm c} {\bm w}_{T,i}^{\rm v} ({\bm \xi}) \leq B, \qquad \text{ and } \qquad \max \;\; \mathbb{E} \; \; \left[  \sum_{t \in \sets T}  \sum_{i \in \sets I} \;\; \theta^{t-1} {\bm \xi}_i^{\rm v} {{\bm z}_{t,i}}({\bm \xi}) \right],
$$
respectively. 
%
%The objective of the BB problem is to maximize the expected profit, i.e.,
%
% $$
% \max \;\; \mathbb{E} \; \; \left[  \sum_{t \in \sets T}  \sum_{i \in \sets I} \;\; \theta^{t-1} {\bm \xi}_i^{\rm v} {{\bm z}_{t,i}}({\bm \xi}) \right].
% $$
%
% The budget constraint can be expressed compactly as
% $$
% \sum_{i \in \sets I} {\bm \xi}_i^{\rm c} {\bm w}_{T,i}^{\rm v} ({\bm \xi}) \leq B.
% $$
% %
% The objective of the BB problem is to maximize the expected profit, i.e.,
% %
% $$
% \max \;\; \mathbb{E} \; \; \left[  \sum_{t \in \sets T}  \sum_{i \in \sets I} \;\; \theta^{t-1} {\bm \xi}_i^{\rm v} {{\bm z}_{t,i}}({\bm \xi}) \right].
% $$
%
The full model for this problem can be found in Electronic Companion~\ref{sec:EC_BB_formulation}.

\subsubsection{Model in RO\cpluspluslogo}
\label{sec:BBcode}
We assume that the data, decision variables, and uncertain parameters of the problem have been defined as in Tables~\ref{tab:BBparams} and~\ref{tab:BBvars}. 
\begin{table}[ht!]
    \centering \renewcommand{\arraystretch}{0.75}
    \begin{tabular}{c|c|c|c}
    \hline
        Model Parameter & \cpluspluslogo{} Variable Name  &  \cpluspluslogo{} Variable Type & \cpluspluslogo{} Map Keys \\
        \hline \hline
         %$T(t)$ & \texttt{T(t)} & \texttt{uint} & NA \\
         %$I(i)$ & \texttt{I(i)} & \texttt{uint} & NA \\
         $B$ & \texttt{B} & \texttt{double} & NA \\
         $\overline{\bm \xi}_i^{\rm c}$,\; $i\in\sets I$ &  \texttt{CostUB} & \texttt{map<uint,double>}& \texttt{i=1\ldots I}  \\
         $\overline{\bm \xi}_i^{\rm v}$,\; $i\in\sets I$ & \texttt{ValueUB} & \texttt{map<uint,double>} & \texttt{i=1\ldots I}  \\
        \hline
    \end{tabular}
    \caption{List of model parameters and their associated \cpluspluslogo{} variables for the BB problem. The parameters $T(t)$ and $I(i)$ are as in Table~\ref{tab:PBparams} and we thus omit them here.}
    \label{tab:BBparams}
\end{table}
%
%\vspace{-0.5cm}
\begin{table}[ht!]
    \centering \renewcommand{\arraystretch}{0.75}
    \begin{tabular}{c|c|c|c}
    \hline
        Parameter & \cpluspluslogo{} Nm.  &  \cpluspluslogo{} Type & \cpluspluslogo{} Map Keys \\
        \hline \hline
         %${\bm z}_{i,t}$, $i\in\sets I$, $t\in\sets T$ &  \texttt{Keep} & \texttt{map<uint,map<uint,ROCPPVarIF\_Ptr> >}& \texttt{1\ldots I, 1\ldots T}  \\
         ${\bm w}_{t,i}^{\rm c}$, $i\in\sets I$, $t\in\sets T$ &  \texttt{MVcost} &
         \texttt{map<uint,map<uint,ROCPPVarIF\_Ptr> >}& \texttt{1\ldots T, 1\ldots I}  \\
         ${\bm w}_{t,i}^{\rm v}$, $i\in\sets I$, $t\in\sets T$ &  \texttt{MVval} &
         \texttt{map<uint,map<uint,ROCPPVarIF\_Ptr> >}& \texttt{1\ldots T, 1\ldots I}  \\
         ${\bm \xi}_i^{\rm c}$, $i\in\sets I$ &  \texttt{Cost} &
         \texttt{map<uint,ROCPPUnc\_Ptr>}& \texttt{i=1\ldots I}  \\
         ${\bm \xi}_i^{\rm v}$, $i\in\sets I$ &  \texttt{Value} &
         \texttt{map<uint,ROCPPUnc\_Ptr>}& \texttt{i=1\ldots I}  \\
        \hline
    \end{tabular}
    \caption{List of model variables and uncertainties and their associated \cpluspluslogo{} variables for the BB problem. The variables ${\bm z}_{i,t}$, $i\in\sets I$, $t\in\sets T$, are as in Table~\ref{tab:PBvars} and we thus omit them here.}
    \label{tab:BBvars}
\end{table}
%\vspace{-0.5cm}

We create a stochastic model with decision-dependent information discovery as follows.
\lstset{firstnumber=1}
\begin{lstlisting}[firstline=1,lastline=2]
// Create an empty stochastic model with T periods for the BB problem
ROCPPOptModelIF_Ptr BBModel(new ROCPPDDUOptModel(T, stochastic));
\end{lstlisting}
To model the requirement that ${\bm \xi}_i^{\rm c}$ and ${\bm \xi}_i^{\rm v}$ must be observed simultaneously, the function \texttt{pair\_uncertainties()} may be employed in \ROCPP.
\lstset{firstnumber=23}
\begin{lstlisting}[firstline=23,lastline=24]
for (uint i = 1; i <= I; i++)
    BBModel->pair_uncertainties(Value[i], Cost[i]);
\end{lstlisting}
To build the budget constraint we use the \ROCPP{} expression \texttt{AmountSpent}.
\lstset{firstnumber=53}
\begin{lstlisting}[firstline=53,lastline=57]
// Constraint on the amount spent
ROCPPExpr_Ptr AmountSpent(new ROCPPExpr());
for (uint i = 1; i <= I; i++)
    AmountSpent = AmountSpent + Cost[i] * MVval[T][i];
BBModel->add_constraint(AmountSpent <= B);
\end{lstlisting}

The construction of the remaining constraints and of the objective parallels that for the Pandora Box problem and we thus omit it. We refer to~\ref{sec:EC_BB_code} for the full code.

\subsubsection{Solution in \ROCPP}
\label{sec:BB_solution}

The BB problem is a multi-stage stochastic problem with decision-dependent information discovery, see~\cite{DDI_VKB}. We thus propose to solve it using PWC decision rules. We consider the instance of BB detailed in Electronic Companion~\ref{sec:EC_BB_instance}, which has $T=4$ and $I=5$. To employ a breakpoint configuration ${\bm r}=(1,1,1,1,1,3,3,1,3,1)$ for the PWC approximation, we use the following code. 
\lstset{firstnumber=74}
\begin{lstlisting}[firstline=74,lastline=80]
// Build the map containing the breakpoint configuration
map<string,uint> BPconfig;
BPconfig["Value_1"] = 3;
BPconfig["Value_2"] = 3;
BPconfig["Value_4"] = 3;
// Construct the PWC decision rule approximator
ROCPPApproximator_Ptr pPWCApprox(new ROCPPPiecewiseApprox(BBModel,BPconfig));
\end{lstlisting}
\ROCPP{} can approximate and robustify the problem in 1 second using lines~\ref{ln:BB_robustify_solve_start}-\ref{ln:BB_robustify_solve_end} in~\ref{sec:EC_BB_code}. Under this breakpoint configuration, the optimal profit is $934.2$, compared to $792.5$ for the static decision rule. The optimal solution can be printed to screen using the \texttt{printOut} function, see lines~\ref{ln:BB_print_start}-\ref{ln:BB_print_end} in~\ref{sec:EC_BB_code}. 
% %
% \lstset{firstnumber=89}
% \lstinputlisting[firstline=89,lastline=92]{bestbox.cc}
% %
Part of the resulting output is
\begin{lstlisting}[numbers=none,language=bash]
On subset 1111131111: Keep_4_1 = 1
Uncertain parameter Value_4 on subset 1111112131 is observed at stage 3
\end{lstlisting}
Thus, on subset ${\bm s}=(1,1,1,1,1,1,3,1,1,1)$, the first box is kept at time~4. On subset ${\bm s}=(1,1,1,1,1,1,2,1,3,1)$, box 4 is opened at time 3 (resp.\ 2).

\section{ROB File Format}
\label{sec:file_format}
Given a robust/stochastic optimization problem expressed in \ROCPP{}, our platform can generate a file displaying the problem in human readable format. We use the Pandora Box problem from Section~\ref{sec:pandora} to illustrate our proposed format, with extension ``.rob''.   

The file starts with the \texttt{Objective} part that presents the objective function of the problem: to minimize either expected or worst-case costs, as indicated by \texttt{E} or \texttt{max}, respectively. For example, since the PB problem optimizes worst-case profit, we obtain the following.
\begin{lstlisting}[numbers=none, language=bash]
Objective:
min max -1 Keep_1_1 Value_1 -1 Keep_1_2 Value_2 -1 Keep_1_3 Value_3  ...
\end{lstlisting}
Then come the \texttt{Constraints} and \texttt{Uncertainty Set} parts, which list the constraints using interpretable ``\lstinline{<=}'', ``\lstinline{>=}'', and ``\lstinline{==}'' operators. We list one constraint for each part here.% for illustration purposes.
\begin{lstlisting}[numbers=none, language=bash]
Constraints:
c0: -1 mValue_2_1 +1 mValue_1_1 <= +0 ...
Uncertainty Set:
c0: -1 Factor_1 <= +1 ...
\end{lstlisting}
The next part, \texttt{Decision Variables}, lists the decision variables of the problem. For each variable, we list its name, type, whether it is static or adaptive, its time stage, and whether it is a measurement variable or not. If it is a measurement variable, we also display the uncertain parameter whose time of revelation it controls.
\begin{lstlisting}[numbers=none, language=bash]
Decision Variables:
Keep_1_1: Boolean, Static, 1, Non-Measurement
mValue_2_2: Boolean, Adaptive, 2, Measurement, Value_2
\end{lstlisting}
The \texttt{Bounds} part then displays the upper and lower bounds for the decision variables.
\begin{lstlisting}[numbers=none, language=bash]
Bounds:
0 <= Keep_1_1 <= 1
\end{lstlisting}
Finally, the \texttt{Uncertainties} part lists, for each uncertain parameter, its name, whether the parameter is observable or not, its time stage, if the parameter has a time of revelation that is decision-dependent, and the first and last stages when the parameter can be observed.
\begin{lstlisting}[numbers=none, language=bash]
Uncertainties:
Factor_4: Not Observable, 1, Non-DDU
Value_1: Observable, 1, DDU, 1, 4
\end{lstlisting}
%%%%%%%%%%%%%%%%%%%%%%%%%%%%%%%%%%%%%%%%%%%%%%%%%%%%%%%%%%%%%%%%%%%%%%%%%%%%%%%%%%%%%%
%%%%%%%%%%%%%%%%%%%%%%%%%%%%%%%%%%%%%%%%%%%%%%%%%%%%%%%%%%%%%%%%%%%%%%%%%%%%%%%%%%%%%%
%%%%%%%%%%%%%%%%%%%%%%%%%%%%%%%%%%%%%% EXTENSIONS %%%%%%%%%%%%%%%%%%%%%%%%%%%%%%%%%%%%
%%%%%%%%%%%%%%%%%%%%%%%%%%%%%%%%%%%%%%%%%%%%%%%%%%%%%%%%%%%%%%%%%%%%%%%%%%%%%%%%%%%%%%
%%%%%%%%%%%%%%%%%%%%%%%%%%%%%%%%%%%%%%%%%%%%%%%%%%%%%%%%%%%%%%%%%%%%%%%%%%%%%%%%%%%%%%

\section{Extensions}
\label{sec:extensions}

\subsection{Integer Decision Variables}
\label{sec:integer_variables}

\ROCPP{} can solve problems involving integer decision variables. In the case of the CDR/PWC approximations, integer adaptive variables are directly approximated by constant/piecewise constant decisions that are integer on each subset. In the case of the finite adaptability approximation, bounded integer variables are automatically expressed as finite sums of binary variables before the finite adaptability approximation is applied.

\subsection{Decision-Dependent Uncertainty Sets}
\label{sec:ddus}

\ROCPP{} can solve problems involving decision-dependent uncertainty sets of the form 
\begin{equation*}
\Xi({\bm z}) \; := \left\{ {\bm \xi} \in \reals^k \; : \; \exists {\bm \zeta}^s \in \reals^{k_s}, \; s=1,\ldots, S \; : \; {\bm P}^s({\bm z}) {\bm \xi} + {\bm Q}^s({\bm z}) {\bm \zeta}^s + {\bm q}^s({\bm z}) \in \sets K^s, \; s=1,\ldots, S \right\},
\end{equation*}
where ${\bm z}$ are static binary variables and ${\bm P}^s({\bm z}) \in \reals^{r_s \times k}$, ${\bm Q}^s({\bm z}) \in \reals^{r_s \times k_s}$, and ${\bm q}^s({\bm z}) \in \reals^{r_s}$, are all linear in ${\bm z}$, and $\sets K^s$, $s=1\ldots,S$, are closed convex pointed cones in $\reals^{r_s}$.

%where ${\bm z}$ are static binary variables and ${\bm V}({\bm z})$ and ${\bm v}({\bm z})$ are linear in ${\bm z}$.

% The problem with decision-dependent uncertainty set~$\Xi := \left\{ {\bm \xi} \; : {\bm A}\left({\bm x}\right){\bm \xi} \leq {\bm b}({\bm x}) \right\}$ is also solvable in \ROCPP. For PWL and LDR approximator, we allow binary and integer variables in the left hand-side expression and all kinds of variables in the right hand-side. For K-Adaptability, only binary variable is allowed.

% \subsection{Approximate $K$-Adaptability for Real-Valued Decisions}
% \label{K-Adapt_real}
% \noteqj{Not sure what to put here}

\subsection{Limited Memory Decision Rules}
\label{sec:limited_memory}

For problems involving long time-horizons ($>100$), the LDR/CDR and PWL/PWC decision rules can become computationally expensive. Limited memory decision rules approximate adaptive decisions by linear functions of the \emph{recent} history of observations. The \texttt{memory} parameter of the \texttt{LDRCDRApproximator} can be used in \ROCPP{} to trade-off optimality with computational complexity.%: the larger (smaller) the memory value, the more (less) flexible the decisions and the greater (lesser) the computational cost. 

%Information may lose after a while, so the decision variable can not depend on the uncertainty observed too distant in the past. To model the above property, our platform provides a limited memory decision rule. Users can decide the number of periods that observed information could be kept by setting the memory parameter in the Piecewise approximator and linear/constant decision rule approximator.

\subsection{Warm Start}
\label{sec:warm_start}

\paragraph{Finite Adaptability.} \ROCPP{} provides the ability to warm start the solution to a finite adaptability problem with $(K_1,\ldots,K_T) \in \integers_{+}^T$ contingency plans using the solution to a finite adaptability problem with fewer contingency plans $(K_1',\ldots,K_T') \in \integers_{+}^T$, where $K_t' \leq K_t$ for all $t \in \sets T$ and $K_t' < K_t$ for at least one $t \in \sets T$.

\paragraph{PWC/PWL Decision Rule.} Similarly, \ROCPP{} provides the ability to warm start the solution to a PWC/PWL decision rule problem with breakpoint configuration ${\bm r} = ({\bm r}_1,\ldots,{\bm r}_k)$ using the solution to a PWC/PWL decision rule problem with breakpoint configuration ${\bm r}' = ({\bm r}_1',\ldots,{\bm r}_k')$. The only requirement is that ${\bm r}_i \geq {\bm r}_i'$ and ${\bm r}_i / {\bm r}_i'$ be a multiple of 2 for all $i \in \{1,\ldots,k\}$ such that ${\bm r}_i > {\bm r}_i'$.

%The warm start method of Gurobi is also available in our platform. Given the solution of a problem with smaller~$K$ or the number of partitions, the K-Adaptability or Piecewise approximator can calculate the warm start map for the same problem but with a larger approximation size. Then we pass the map into the Gurobi interface to set the \texttt{GRB\_DoubleAttr\_Start} to construct an initial solution for each variable.

\newpage

% \newpage
% \begingroup
% \parindent 0pt
% \parskip 2ex
\def\enotesize{\small}
\theendnotes
%\endgroup

% Acknowledgments here
\ACKNOWLEDGMENT{%
This work was supported in part by the Operations Engineering Program of the National Science Foundation under NSF Award No.\ 1763108. The authors are grateful to Daniel Kuhn, Ber{\c{c}} Rustem, and Wolfram Wiesemann, for valuable discussions that helped shape this work.
}% 

%%%%%%%%%%%%%%%%%%%%%%%%%%%%%%%%%%%%%%%%%%%%%%%%%%%%%%%%%%%%%%%%%%%%%%%%%%%%%%%%%%%%%
%%%%%%%%%%%%%%%%%%%%%%%%%%%%%%%%%%%%%%%%%%%%%%%%%%%%%%%%%%%%%%%%%%%%%%%%%%%%%%%%%%%%%

%%%%%%%%%%%%%%%%%%%%%%%%%%%%%%%%%%%%%%%%%%%%%%%%%%%%%%%%%%%%%%%%%%%%%%%%%%%%%%%%%%%%%
%%%%%%%%%%%%%%%%%%%%%%%%%%%%%%%%%%%%%%%%%%%%%%%%%%%%%%%%%%%%%%%%%%%%%%%%%%%%%%%%%%%%%
%%%%%%%%%%%%%%%%%%%%%%%%%%%%%%%%%%%%%%%% BIBLIOGRAPHY %%%%%%%%%%%%%%%%%%%%%%%%%%%%%%%
%%%%%%%%%%%%%%%%%%%%%%%%%%%%%%%%%%%%%%%%%%%%%%%%%%%%%%%%%%%%%%%%%%%%%%%%%%%%%%%%%%%%%
%%%%%%%%%%%%%%%%%%%%%%%%%%%%%%%%%%%%%%%%%%%%%%%%%%%%%%%%%%%%%%%%%%%%%%%%%%%%%%%%%%%%%

\bibliographystyle{informs2014} 
%\bibliographystyle{authoryear}
%\bibliography{Mendeley.bib}
%\bibliography{bibliography_noURL.bib}

%%%%%%%%%%%%%%%%%%%%%%%%%%%%%%%%%%%%%%%%%%%%%%%%%%%%%%%%%%%%%%%%%%%%%%%%%%%%%%%%%%%%%
%%%%%%%%%%%%%%%%%%%%%%%%%%%%%%%%%%%%%%%%%%%%%%%%%%%%%%%%%%%%%%%%%%%%%%%%%%%%%%%%%%%%%
%%%%%%%%%%%%%%%%%%%%%%%%%%%%%%%%%%%%%%%% E-COMPANION %%%%%%%%%%%%%%%%%%%%%%%%%%%%%%%%
%%%%%%%%%%%%%%%%%%%%%%%%%%%%%%%%%%%%%%%%%%%%%%%%%%%%%%%%%%%%%%%%%%%%%%%%%%%%%%%%%%%%%
%%%%%%%%%%%%%%%%%%%%%%%%%%%%%%%%%%%%%%%%%%%%%%%%%%%%%%%%%%%%%%%%%%%%%%%%%%%%%%%%%%%%%

% Here starts the e-companion (EC)
%%%%%%%%%%%%%%%%%%%%%%%%%%%%%%%%%%%%%%%%%%%%%%%%%%%%%%%%%
\ECSwitch

%\ECDisclaimer
%%%%%%%%%%%%%%%%%%%%%%%%%%%%%%%%%%%%%%%%%%%%%%%%%%%%%%%%%

%% Main head for the e-companion
\ECHead{E-Companion}

%%%%%%%%%%%%%%%%%%%%%%%%%%%%%%%%%%%%%%%%%%%%%%%%%%%%%%%%%
%%%%%%%%%%%%%%%%%%%%%% RSFC %%%%%%%%%%%%%%%%%%%%%%%%%%%%%
%%%%%%%%%%%%%%%%%%%%%%%%%%%%%%%%%%%%%%%%%%%%%%%%%%%%%%%%%

\section{Supplemental Material: Retailer-Supplier Problem}
\label{sec:EC_RSFC}

%%%%%%%%%%%%%%%%%%%%%%%%%%%%%%%%%%%%%%%%%%%%%%%%%%%%%%%%%

\subsection{Retailer-Supplier Problem: Mathematical Formulation}
\label{sec:EC_RSFC_formulation}

Using the notation introduced in Section~\ref{sec:retailer}, the robust RSFC problem can be expressed mathematically as:
\begin{equation*}\renewcommand{\arraystretch}{\pveqnstretch}
    \begin{array}{cl}
        \minimize & \quad \displaystyle \max_{{\bm \xi} \in \Xi} \; \sum_{t \in \sets T} {\bm c}_t^{\rm o} {\bm y}_t^{\rm o}({\bm \xi}) +  {\bm y}_t^{\rm dc}({\bm \xi}) + {\bm y}_t^{\rm dp}({\bm \xi}) + {\bm y}_{t+1}^{\rm hs}({\bm \xi})  \\
        \subjectto & \quad {\bm y}_t^{\rm c} \in \reals_+ \quad \forall t \in \sets T \\
        & \quad  {\bm y}_t^{\rm o} , \; {\bm y}_t^{\rm dc}, \; {\bm y}_t^{\rm dp}, \; {\bm y}_{t+1}^{\rm hs} \in \mathcal L_T^1  \quad \forall t \in \sets T \\
        & \quad \left. \!\!\!\begin{array}{l}
        {\bm x}^{\rm i}_{t+1}({\bm \xi}) \; = \; {\bm x}^{\rm i}_t({\bm \xi}) + {\bm y}^{\rm o}_t({\bm \xi}) - {\bm \xi}_{t+1} \\
        \underline{\bm y}_t^{\rm o} \leq {\bm y}_t^{\rm o}({\bm \xi}) \leq \overline{\bm y}_t^{\rm o} , \;\; \underline{\bm y}_t^{\rm co} \leq \sum_{\tau = 1}^t {\bm y}_\tau^{\rm o}({\bm \xi}) \leq \overline{\bm y}_t^{\rm co} \\
        {\bm y}_t^{\rm dc} \geq {\bm c}_t^{\rm dc+} ( {\bm y}_t^{\rm c} - {\bm y}_{t-1}^{\rm c} ) \\
        {\bm y}_t^{\rm dc} \geq {\bm c}_t^{\rm dc-} ({\bm y}_{t-1}^{\rm c} - {\bm y}_{t}^{\rm c} ) \\
        {\bm y}_t^{\rm dp} \geq {\bm c}_t^{\rm dp+} ({\bm y}_t^{\rm o}({\bm \xi}) - {\bm y}_{t}^{\rm c}) \\
        {\bm y}_t^{\rm dp}({\bm \xi}) \geq {\bm c}_t^{\rm dp-} ({\bm y}_t^{\rm c} - {\bm y}_{t}^{\rm o}({\bm \xi})) \\
        {\bm y}_{t+1}^{\rm hs}({\bm \xi}) \geq {\bm c}_{t+1}^{\rm h}{\bm x}^{\rm i}_{t+1}({\bm \xi})   \\
        {\bm y}_{t+1}^{\rm hs}({\bm \xi}) \geq -{\bm c}_{t+1}^{\rm s} {\bm x}^{\rm i}_{t+1}({\bm \xi})
        \end{array} \quad \right\} \quad \forall t \in \sets T , \; {\bm \xi} \in \Xi \\
        & \quad \left. \!\!\! \begin{array}{l}
        {\bm y}_t^{\rm o}({\bm \xi}) = {\bm y}_t^{\rm o}({\bm \xi}') \\
        {\bm y}_t^{\rm dc}({\bm \xi}) = {\bm y}_t^{\rm dc}({\bm \xi}') \\
        {\bm y}_t^{\rm dp}({\bm \xi}) = {\bm y}_t^{\rm dp}({\bm \xi}')
        \end{array} \qquad \quad \right\} \quad \forall {\bm \xi}, {\bm \xi}' \in \Xi : {\bm w}_{t-1} \circ {\bm \xi} = {\bm w}_{t-1} \circ {\bm \xi}', \; \forall t \in \sets T \\
        & \quad  {\bm y}_{t+1}^{\rm hs}({\bm \xi}) = {\bm y}_{t+1}^{\rm hs}({\bm \xi}') \qquad \quad \forall {\bm \xi}, {\bm \xi}' \in \Xi : {\bm w}_{t} \circ {\bm \xi} = {\bm w}_{t} \circ {\bm \xi}', \; \forall t \in \sets T.
    \end{array} 
\end{equation*}
The last set of constraints corresponds to non-anticipativity constraints.  The other constraints are explained in Section~\ref{sec:RSFC_description}.

%%%%%%%%%%%%%%%%%%%%%%%%%%%%%%%%%%%%%%%%%%%%%%%%%%%%%%%%%

\subsection{Retailer-Supplier Problem: Full \ROCPP{} Code}
\label{sec:EC_RSFC_code}

\lstset{firstnumber=1}
\begin{lstlisting}
// Create an empty robust model with T + 1 periods for the RSFC problem
ROCPPOptModelIF_Ptr RSFCModel(new ROCPPUncOptModel(T+1, robust));
// Create the Demand map to store the uncertain parameters of the problem
map<uint,ROCPPUnc_Ptr> Demand;
// Iterate over all time periods when there is uncertainty
for (uint t=1; t<=T; t++)
    // Create the uncertain parameters, and add them to Demand
    Demand[t] = ROCPPUnc_Ptr(new ROCPPUnc("Demand_"+to_string(t),t+1));
// Create maps to store the decision variables of the problem
map<uint,ROCPPVarIF_Ptr> Orders, Commits;    // Quantity ordered, Commitments made
// Iterate over all time periods from 1 to T
for (uint t=1; t<=T; t++) {
    // Create the commitment variables (these are static)
    Commits[t]=ROCPPVarIF_Ptr(new ROCPPStaticVarReal("Commit_"+to_string(t),0.));
    if (t==1) // In the first period, the order variables are static
        Orders[t] = ROCPPVarIF_Ptr(new ROCPPStaticVarReal("Order_"+to_string(t),OrderLB[t],OrderUB[t]));
    else // In the other periods, the order variables are adaptive
        Orders[t] = ROCPPVarIF_Ptr(new ROCPPAdaptVarReal("Order_"+to_string(t),t,OrderLB[t],OrderUB[t]));
}
map<uint,ROCPPVarIF_Ptr> MaxDC; |\label{ln:retailer_vars_2_start}| // Upper bound on deviation between commitments  
map<uint,ROCPPVarIF_Ptr> MaxDP; // Upper bound on deviation from plan
map<uint,ROCPPVarIF_Ptr> MaxHS; // Upper bound on holding and shortage costs
// Iterate over all time periods 1 to T
for (uint t=1; t<=T; t++) {
    // Create upper bounds on the deviation between successive commitments
    MaxDC[t] = ROCPPVarIF_Ptr(new ROCPPStaticVarReal("MaxDC_"+to_string(t)));
    // Create upper bounds on the deviation of orders from commitments
    if (t==1)   // In the first period, these are static
        MaxDP[t] = ROCPPVarIF_Ptr(new ROCPPStaticVarReal("MaxDP_"+to_string(t)));
    else        // In the other periods, these are adaptive
        MaxDP[t] = ROCPPVarIF_Ptr(new ROCPPAdaptVarReal("MaxDP_"+to_string(t),t));
    // Create upper bounds on holding and shortage costs (these are adaptive)
    MaxHS[t+1]=ROCPPVarIF_Ptr(new ROCPPAdaptVarReal("MaxHS_"+to_string(t+1),t+1));
} |\label{ln:retailer_vars_2_end}|
// Create the constraints of the problem
// Create an expression for the amount of inventory held and initialize it
ROCPPExpr_Ptr Inventory(new ROCPPExpr());
Inventory = Inventory + InitInventory;
// Create an expression for the cumulative amount ordered
ROCPPExpr_Ptr CumOrders(new ROCPPExpr());
// Iterate over all time periods and add the constraints to the problem
for (uint t=1; t<=T; t++) {
    // Create the upper and lower bounds on the cumulative orders
    CumOrders = CumOrders + Orders[t];
    RSFCModel->add_constraint(CumOrders >= CumOrderLB[t]);
    RSFCModel->add_constraint(CumOrders <= CumOrderUB[t]);
    // Update the inventory
    Inventory = Inventory + Orders[t] - Demand[t];
    // Create upper bound on shortage/holding costs
    RSFCModel->add_constraint(MaxHS[t+1] >= HoldingCost[t+1]*Inventory);
    RSFCModel->add_constraint(MaxHS[t+1] >= (-1.*ShortageCost[t+1]*Inventory));
}
|\label{ln:retailer_extra_cstr_start}| // Iterate over all time periods and add the constraints to the problem
for (uint t=1; t<=T; t++) {
    // Create upper bound on deviations from commitments
    RSFCModel->add_constraint(MaxDP[t] >= CostDPp*(Orders[t]-Commits[t]));
    RSFCModel->add_constraint(MaxDP[t] >= -CostDPm*(Orders[t]-Commits[t]));
    // Create upper bound on deviations between commitments
    if (t==1) {
        RSFCModel->add_constraint(MaxDC[t] >= CostDCp*(Commits[t]-InitCommit));
        RSFCModel->add_constraint(MaxDC[t] >= -CostDCm*(Commits[t]-InitCommit));
    }
    else {
        RSFCModel->add_constraint(MaxDC[t] >= CostDCp*(Commits[t]-Commits[t-1]));
        RSFCModel->add_constraint(MaxDC[t] >= -CostDCm*(Commits[t]-Commits[t-1]));
    }
} |\label{ln:retailer_extra_cstr_end}|
// Create an expression that will contain the objective function
ROCPPExpr_Ptr RSFCObj(new ROCPPExpr());
// Iterate over all periods and add the terms to the objective function
for (uint t=1; t<=T; t++)
    RSFCObj = RSFCObj + OrderCost*Orders[t] + MaxDC[t] + MaxDP[t] + MaxHS[t+1];
RSFCModel->set_objective(RSFCObj); // Add the objective to the problem
for (uint t=1; t<=T; t++) |\label{ln:uset_start}|{
    // Add the upper and lower bounds on the demand to the uncertainty set
    RSFCModel->add_constraint_uncset(Demand[t] >= NomDemand*(1.0-rho));
    RSFCModel->add_constraint_uncset(Demand[t] <= NomDemand*(1.0+rho));
} |\label{ln:uset_end}|
// Construct the linear/constant decision rule approximator
ROCPPApproximator_Ptr pLDRApprox(new ROCPPLCDRApprox(RSFCModel));
// Approximate the adaptive decisions using the LDR/CDR approximator 
ROCPPMISOCP_Ptr RSFCModelLDR(pLDRApprox->DoMyThing(RSFCModel));
// Construct the solver (in this case, use gurobi as deterministic solver)
ROCPPSolver_Ptr pSolver(new ROCPPGurobi(SolverParams()));
// Solve the problem
pSolver->solve(RSFCModelLDR);
// Retrieve the optimal solution from the solver
map<string,double> optimalSln(pSolver->getSolution());
// Print the optimal decision (from the original model)
pLDRApprox->printOut(RSFCModelLDR, optimalSln, Orders[10]));
// Get the optimal objective value
double optVal(pSolver->getOptValue());
\end{lstlisting}

%%%%%%%%%%%%%%%%%%%%%%%%%%%%%%%%%%%%%%%%%%%%%%%%%%%%%%%%%

\subsection{Retailer-Supplier Problem: Instance Parameters}
\label{sec:EC_RSFC_instance}

The parameters for the instance of the problem that we solve in Section~\ref{sec:RSFC_solution} are provided in Table~\ref{tab:RSFC_instance_params}. They correspond to the data from instance W12 in~\cite{BenTal_RSFC}.
\begin{table}[ht!]
    \centering
    \begin{tabular}{|c|c|c|c|c|c|c|c|c|c|c|c|c|c|c|c|c|c}
        \hline
        $T$ & ${\bm x}^{\rm i}_1$ & ${\bm y}^{\rm c}_0$ & $\overline \xi$ & $\rho$ & $\underline{\bm y}_t^{\rm o}$ & $\overline{\bm y}_t^{\rm o}$ & $\underline{\bm y}_t^{\rm co}$ & $\overline{\bm y}_t^{\rm co}$ & ${\bm c}^{\rm o}_{t}$ & ${\bm c}^{\rm h}_{t+1}$ & ${\bm c}^{\rm s}_{t+1}$ & ${\bm c}^{\rm dc+}_{t}$ & ${\bm c}^{\rm dc-}_{t}$ & ${\bm c}^{\rm dp+}_{t}$ & ${\bm c}^{\rm dp-}_{t}$ \\
        \hline \hline
        12 & 0 & 100 & 100 & 10\% & 0 & 200 & 0 & $200t$ & 10 & 2 & 10 & 10 & 10 & 10 & 10 \\
        \hline
    \end{tabular}
    \caption{Parameters for the instance of the RSFS problem that we solve in Section~\ref{sec:RSFC_solution}.}
    \label{tab:RSFC_instance_params}
\end{table}

%%%%%%%%%%%%%%%%%%%%%%%%%%%%%%%%%%%%%%%%%%%%%%%%%%%%%%%%%
%%%%%%%%%%%%%%%%%%%%%%%% PB %%%%%%%%%%%%%%%%%%%%%%%%%%%%%
%%%%%%%%%%%%%%%%%%%%%%%%%%%%%%%%%%%%%%%%%%%%%%%%%%%%%%%%%

\section{Supplemental Material: Robust Pandora's Box Problem}
\label{sec:EC_PB}

\subsection{Robust Pandora's Box Problem: Mathematical Formulation}
\label{sec:EC_PB_formulation}

Using the notation introduced in Section~\ref{sec:pandora}, the robust PB problem can be expressed mathematically as:
\begin{equation*}\renewcommand{\arraystretch}{\pveqnstretch}
    \begin{array}{cl}
        \maximize & \quad  \displaystyle \min_{{\bm \xi} \in \Xi} \;\; \sum_{t \in \sets T} \sum_{i \in \sets I} \theta^{t-1} {\bm \xi}_i {\bm z}_{t,i}({\bm \xi}) - {\bm c}_i ({\bm w}_{t,i}({\bm \xi})-{\bm w}_{t-1,i}({\bm \xi}))  \\
        \subjectto & \quad  {\bm z}_{t,i}, \; {\bm w}_{t,i} \in \{0,1\} \quad \forall t \in \sets T, \; \forall i \in \sets I  \\
        & \quad  \left. \!\!\!\begin{array}{l}
        \displaystyle \sum_{i \in \sets I} ( {\bm w}_{t,i}({\bm \xi}) - {\bm w}_{t-1,i}({\bm \xi}) ) \; \leq \; 1 - \sum_{\tau=1}^t {\bm z}_{t,i}({\bm \xi}) \\
        {\bm z}_{t,i}({\bm \xi}) \; \leq \; {\bm w}_{t-1,i}({\bm \xi}) \quad \forall i \in \sets I\\
        \end{array} \quad \right\} \quad \forall t \in \sets T , \; {\bm \xi} \in \Xi \\
        & \quad  \left. \!\!\! \begin{array}{l}
        {\bm z}_{t,i}({\bm \xi}) = {\bm z}_{t,i}({\bm \xi}') \\
        {\bm w}_{t,i}({\bm \xi}) = {\bm w}_{t,i}({\bm \xi}') \\
        \end{array} \qquad \quad \right\} \quad \forall {\bm \xi}, {\bm \xi}' \in \Xi : {\bm w}_{t-1} \circ {\bm \xi} = {\bm w}_{t-1} \circ {\bm \xi}', \; \forall i \in \sets I, \; \forall t \in \sets T .
    \end{array} 
\end{equation*}
The last set of constraints in this problem are non-anticipativity constraints. The other constraints are explained in Section~\ref{sec:RPBdescription}.

%%%%%%%%%%%%%%%%%%%%%%%%%%%%%%%%%%%%%%%%%%%%%%%%%%%%%%%%%

\subsection{Robust Pandora's Box Problem: Full \ROCPP{} Code}
\label{sec:EC_PB_code}

\lstset{firstnumber=1}
\begin{lstlisting}
// Create an empty robust model with T periods for the PB problem
ROCPPOptModelIF_Ptr PBModel(new ROCPPDDUOptModel(T, robust));
// Create empty maps to store the uncertain parameters
map<uint, ROCPPUnc_Ptr> Value, Factor;
for (uint i = 1; i <= I; i++)
    // Create the uncertainty associated with box i and add it to Value
    Value[i] = ROCPPUnc_Ptr(new ROCPPUnc("Value_"+to_string(i)));
for (uint m = 1; m <= M; m++)
    // The risk factors are not observable
    Factor[m]= ROCPPUnc_Ptr(new ROCPPUnc("Factor_"+to_string(m),1,false));
map<uint, map<uint, ROCPPVarIF_Ptr> > MeasVar;
for (uint i = 1; i <= I; i++) {
    // Create the measurement variables associated with the value of box i
    PBModel->add_ddu(Value[i], 1, T, obsCost[i]);
    // Get the measurement variables and store them in MeasVar
    for (uint t = 1; t <= T; t++)
        MeasVar[t][i] = PBModel->getMeasVar(Value[i]->getName(), t);
}
map<uint, map<uint, ROCPPVarIF_Ptr> > Keep; |\label{ln:keep_vars_start}|
for (uint t = 1; t <= T; t++) {
    for (uint i = 1; i <= I; i++) {
        if (t == 1)  // In the first period, the Keep variables are static
            Keep[t][i] = ROCPPVarIF_Ptr(new ROCPPStaticVarBool("Keep_"+to_string(t)+"_"+to_string(i)));
        else   // In the other periods, the Keep variables are adaptive
            Keep[t][i] = ROCPPVarIF_Ptr(new ROCPPAdaptVarBool("Keep_"+to_string(t)+"_"+to_string(i), t));
    }
} |\label{ln:keep_vars_end}|
|\label{ln:pandora_cstr_start}| // Create the constraints and add them to the problem
ROCPPExpr_Ptr StoppedSearch(new ROCPPExpr());
for (uint t = 1; t <= T; t++) {
    // Create the constraint that at most one box be opened at t (none if the search has stopped)
    ROCPPExpr_Ptr NumOpened(new ROCPPExpr());
    // Update the expressions and and the constraint to the problem
    for (uint i = 1; i <= I; i++) {
        StoppedSearch = StoppedSearch + Keep[t][i];
        if (t>1)
            NumOpened = NumOpened + MeasVar[t][i] - MeasVar[t-1][i];
        else
            NumOpened = NumOpened + MeasVar[t][i];
    }
    PBModel->add_constraint( NumOpened <= 1. - StoppedSearch );
    // Constraint that only one of the open boxes can be kept
    for (uint i = 1; i <= I; i++)
        PBModel->add_constraint( (t>1) ? (Keep[t][i] <= MeasVar[t-1][i]) : (Keep[t][i] <= 0.));
} |\label{ln:pandora_cstr_end}|
// Create the uncertainty set constraints and add them to the problem
// Add the upper and lower bounds on the risk factors
for (uint m = 1; m <= M; m++) {
    PBModel->add_constraint_uncset(Factor[m] >= -1.0);
    PBModel->add_constraint_uncset(Factor[m] <= 1.0);
}
// Add the expressions for the box values in terms of the risk factors
for (uint i = 1; i <= I; i++) {
    ROCPPExpr_Ptr ValueExpr(new ROCPPExpr());
    for (uint m = 1; m <= M; m++)
        ValueExpr = ValueExpr + RiskCoeff[i][m]*Factor[m];
    PBModel->add_constraint_uncset(  Value[i] == (1.+0.5*ValueExpr) * NomVal[i] );
}
// Create the objective function expression
ROCPPExpr_Ptr PBObj(new ROCPPExpr()); 
for (uint t = 1; t <= T; t++)
    for (uint i = 1; i <= I; i++)
        PBObj = PBObj + pow(theta,t-1)*Value[i]*Keep[t][i];
// Set objective (multiply by -1 for maximization)
PBModel->set_objective(-1.0*PBObj);
// Construct the finite adaptability approximator
ROCPPApproximator_Ptr pKAdaptaApprox(new ROCPPKadaptApprox(PBModel,K));
// Approximate the decisions using finite adaptability and robustify
ROCPPMISOCP_Ptr PBModelKAdapt(pKAdaptaApprox->DoMyThing(PBModel));
// Construct the solver and solve the problem
ROCPPSolver_Ptr pSolver(new ROCPPGurobi(SolverParams()));
pSolver->solve(PBModelKAdapt);
// Retrieve the optimal solution from the solver
map<string,double> optimalSln(pSolver->getSolution());
// Print the optimal decision (from the original model)
pKAdaptaApprox->printOut(PBModelKAdapt, optimalSln, Keep[4][2]);
pKAdaptaApprox->printOut(PBModel, optimalSln, Value[4]);
\end{lstlisting}

%%%%%%%%%%%%%%%%%%%%%%%%%%%%%%%%%%%%%%%%%%%%%%%%%%%%%%%%%

\subsection{Robust Pandora's Box Problem: Instance Parameters}
\label{sec:EC_PB_instance}

The parameters for the instance of the robust Pandora's box problem that we solve in Section~\ref{sec:pandora_solution} are provided in Table~\ref{tab:PB_instance_params}.
\begin{table}[ht!]
    \centering
    \begin{tabular}{|c|c|}
    \hline
    Parameter & Value \\
        \hline  \hline
        $(\theta,T, I,M)$ & (1,4,5,4) \\
        ${\bm c}$ &  $(0.69, 0.43, 0.01, 0.91, 0.64)$ \\ 
        $\overline{\bm \xi}$ & $(5.2, 8, 19.4, 9.6, 13.2)$ \\
        ${\bm \Phi}$ & $\begin{pmatrix}
        0.17 & -0.7 & -0.13 & -0.6 \\
        0.39 & 0.88 & 0.74 & 0.78 \\
        0.17 & -0.6 & -0.17 & -0.84 \\
        0.09 & -0.07 & -0.52 & 0.88 \\
        0.78 & 0.94 & 0.43 & -0.58
        \end{pmatrix}$ \\
        \hline
    \end{tabular}
    \caption{Parameters for the instance of the PB problem that we solve in Section~\ref{sec:pandora_solution}.}
    \label{tab:PB_instance_params}
\end{table}

%%%%%%%%%%%%%%%%%%%%%%%%%%%%%%%%%%%%%%%%%%%%%%%%%%%%%%%%%
%%%%%%%%%%%%%%%%%%%%%%%% BB %%%%%%%%%%%%%%%%%%%%%%%%%%%%%
%%%%%%%%%%%%%%%%%%%%%%%%%%%%%%%%%%%%%%%%%%%%%%%%%%%%%%%%%

\section{Supplemental Material: Stochastic Best Box Problem}
\label{sec:EC_BB}

\subsection{Stochastic Best Box: Problem Formulation}
\label{sec:EC_BB_formulation}

Using the notation introduced in Section~\ref{sec:best_box}, the BB problem can be expressed mathematically as:
\begin{equation*}\renewcommand{\arraystretch}{\pveqnstretch}
    \begin{array}{cl}
        \maximize & \quad  \displaystyle \mathbb E \left[  \sum_{t \in \sets T} \sum_{i \in \sets I} \theta^{t-1} {\bm \xi}_i^{\rm v} {\bm z}_{t,i}({\bm \xi}) \right]  \\
        \subjectto & \quad  {\bm z}_{t,i}, \; {\bm w}_{t,i}^{\rm c}, \; {\bm w}_{t,i}^{\rm v} \in \{0,1\} \quad \forall t \in \sets T, \; \forall i \in \sets I  \\
        & \quad {\bm w}_{t,i}^{\rm c}({\bm \xi}) \; = \; {\bm w}_{t,i}^{\rm v}({\bm \xi}) \quad \forall t \in \sets T, \; \forall i \in \sets I  \\
        & \quad  \left. \!\!\!\begin{array}{l}
        \displaystyle \sum_{i \in \sets I} {\bm \xi}_i^{\rm c} {\bm w}_{T,i}^{\rm v}({\bm \xi}) \leq B \\
        \displaystyle \sum_{i \in \sets I} ( {\bm w}_{t,i}^{\rm v}({\bm \xi}) - {\bm w}_{t-1,i}^{\rm v}({\bm \xi}) ) \; \leq \; 1 - \sum_{\tau=1}^t {\bm z}_{t,i}({\bm \xi}) \\
        {\bm z}_{t,i}({\bm \xi}) \; \leq \; {\bm w}_{t-1,i}^{\rm v}({\bm \xi}) \quad \forall i \in \sets I\\
        \end{array} \quad \right\} \quad \forall t \in \sets T , \; {\bm \xi} \in \Xi \\
        & \quad  \left. \!\!\! \begin{array}{l}
        {\bm z}_{t,i}({\bm \xi}) = {\bm z}_{t,i}({\bm \xi}') \\
        {\bm w}_{t,i}^{\rm c}({\bm \xi}) = {\bm w}_{t,i}^{\rm c}({\bm \xi}') \\
        {\bm w}_{t,i}^{\rm v}({\bm \xi}) = {\bm w}_{t,i}^{\rm v}({\bm \xi}') 
        \end{array} \qquad \quad \right\} \quad \forall {\bm \xi}, {\bm \xi}' \in \Xi : {\bm w}_{t-1} \circ {\bm \xi} = {\bm w}_{t-1} \circ {\bm \xi}', \; \forall i \in \sets I, \; \forall t \in \sets T .
    \end{array} 
\end{equation*}
The first set of constraints stipulates that ${\bm \xi}_i^{\rm c}$ and ${\bm \xi}_i^{\rm v}$ must be observed simultaneously. The second set of constraints is the budget constraint. The third set of constraints stipulates that at each stage, we can either open a box or stop the search, in which case we cannot open a box in the future. The fourth set of constraints ensures that we can only keep a box that we have opened. The last set of constraints correspond to decision-dependent non-anticipativity constraints.

%%%%%%%%%%%%%%%%%%%%%%%%%%%%%%%%%%%%%%%%%%%%%%%%%%%%%%%%%

\subsection{Stochastic Best Box Problem: Full \ROCPP{} Code}
\label{sec:EC_BB_code}

\lstset{firstnumber=1}
\begin{lstlisting}
// Create an empty stochastic model with T periods for the BB problem
ROCPPOptModelIF_Ptr BBModel(new ROCPPDDUOptModel(T, stochastic));
map<uint, ROCPPUnc_Ptr> Value, Cost;
for (uint i = 1; i <= I; i++) {
    // Create the value and cost uncertainties associated with box i
    Value[i] = ROCPPUnc_Ptr(new ROCPPUnc("Value_"+to_string(i)));
    Cost[i] = ROCPPUnc_Ptr(new ROCPPUnc("Cost_"+to_string(i)));
}
// Create the measurement decisions and pair the uncertain parameters
map<uint, map<uint, ROCPPVarIF_Ptr> > MVcost, MVval;
for (uint i = 1; i <= I; i++) {
    // Create the measurement variables associated with the value of box i
    BBModel->add_ddu(Value[i], 1, T, obsCost);
    // Create the measurement variables associated with the cost of box i
    BBModel->add_ddu(Cost[i], 1, T, obsCost);
    // Get the measurement variables and store them in MVval and MVcost
    for (uint t = 1; t <= T; t++) {
        MVval[t][i] = BBModel->getMeasVar(Value[i]->getName(), t);
        MVcost[t][i] = BBModel->getMeasVar(Cost[i]->getName(), t);
    }
}
// Pair the uncertain parameters to ensure they are observed at the same time
for (uint i = 1; i <= I; i++)
    BBModel->pair_uncertainties(Value[i], Cost[i]);
// Create the keep decisions
map<uint, map<uint, ROCPPVarIF_Ptr> > Keep;
for (uint t = 1; t <= T; t++) {
    for (uint i = 1; i <= I; i++) {
        if (t == 1)  // In the first period, the Keep variables are static
            Keep[t][i] = ROCPPVarIF_Ptr(new ROCPPStaticVarBool("Keep_"+to_string(t)+"_"+to_string(i)));
        else   // In the other periods, the Keep variables are adaptive
            Keep[t][i] = ROCPPVarIF_Ptr(new ROCPPAdaptVarBool("Keep_"+to_string(t)+"_"+to_string(i), t));
    }
}
// Create the constraints and add them to the problem
ROCPPExpr_Ptr StoppedSearch(new ROCPPExpr());
for (uint t = 1; t <= T; t++) {
    // Create the constraint that at most one box be opened at t (none if the search has stopped)
    ROCPPExpr_Ptr NumOpened(new ROCPPExpr());
    // Update the expressions and and the constraint to the problem
    for (uint i = 1; i <= I; i++) {
        StoppedSearch = StoppedSearch + Keep[t][i];
        if (t>1)
            NumOpened = NumOpened + MVval[t][i] - MVval[t-1][i];
        else
            NumOpened = NumOpened + MVval[t][i];
    }
    BBModel->add_constraint( NumOpened <= 1. - StoppedSearch );
    // Constraint that only one of the open boxes can be kept
    for (uint i = 1; i <= I; i++)
        BBModel->add_constraint( (t>1) ? (Keep[t][i] <= MVval[t-1][i]) : (Keep[t][i] <= 0.));
}
// Constraint on the amount spent
ROCPPExpr_Ptr AmountSpent(new ROCPPExpr());
for (uint i = 1; i <= I; i++)
    AmountSpent = AmountSpent + Cost[i] * MVval[T][i];
BBModel->add_constraint(AmountSpent <= B);
// Create the uncertainty set constraints and add them to the problem
for (uint i = 1; i <= I; i++) {
    // Add the upper and lower bounds on the values
    BBModel->add_constraint_uncset(Value[i] >= 0.);
    BBModel->add_constraint_uncset(Value[i] <= ValueUB[i]);
    // Add the upper and lower bounds on the costs
    BBModel->add_constraint_uncset(Cost[i] >= 0.);
    BBModel->add_constraint_uncset(Cost[i] <= CostUB[i]);
}
// Create the objective function expression
ROCPPExpr_Ptr BBObj(new ROCPPExpr());
for (uint t = 1; t <= T; t++)
    for (uint i = 1; i <= I; i++)
        BBObj = BBObj + pow(theta,t-1)*Value[i]*Keep[t][i];
// Set objective (multiply by -1 for maximization)
BBModel->set_objective(-1.0*BBObj);
// Build the map containing the breakpoint configuration
map<string,uint> BPconfig;
BPconfig["Value_1"] = 3;
BPconfig["Value_2"] = 3;
BPconfig["Value_4"] = 3;
// Construct the PWC decision rule approximator
ROCPPApproximator_Ptr pPWCApprox(new ROCPPPiecewiseApprox(BBModel,BPconfig));
|\label{ln:BB_robustify_solve_start}| // Approximate the decisions using PWC decision rules and robustify
ROCPPMISOCP_Ptr BBModelPWC(pPWCApprox->DoMyThing(BBModel));
// Construct the solver; in this case, use the gurobi solver as a deterministic solver
ROCPPSolver_Ptr pSolver(new ROCPPGurobi(SolverParams()));
// Solve the problem
pSolver->solve(BBModelPWC); |\label{ln:BB_robustify_solve_end}|
|\label{ln:BB_print_start}| // Retrieve the optimal solution from the solver
map<string,double> optimalSln(pSolver->getSolution());
// Print the optimal decision (from the original model)
// Prints decision rules from the original problem automatically
pPWCApprox->printOut(BBModelPWC, optimalSln, Keep[4][1]);
pPWCApprox->printOut(BBModel, optimalSln, Value[4]); |\label{ln:BB_print_end}|
\end{lstlisting}

%%%%%%%%%%%%%%%%%%%%%%%%%%%%%%%%%%%%%%%%%%%%%%%%%%%%%%%%%

\subsection{Stochastic Best Box: Instance Parameters}
\label{sec:EC_BB_instance}

The parameters for the instance of the stochastic best box problem that we solve in Section~\ref{sec:BB_solution} are provided in Table~\ref{tab:BB_instance_params}.
\begin{table}[ht!]
    \centering
    \begin{tabular}{|c|c|c|c|c|c|}
        \hline
        $T$ & $I$ & $B$ & $\theta$ & $\overline{\bm \xi}^{\rm c}$ & $\overline{\bm \xi}^{\rm v}$ \\
        \hline \hline
        4 & 5 & 163 & 1 & (40,86,55,37,30) & (1030,1585,971,971,694) \\
        \hline
    \end{tabular}
    \caption{Parameters for the instance of the BB problem that we solve in Section~\ref{sec:BB_solution}.}
    \label{tab:BB_instance_params}
\end{table}

\end{document}